\documentclass[12pt]{article}
\usepackage{amssymb,amscd,amsmath,amsthm,fullpage}
\makeatletter
 
\def\diagram{\m@th\leftwidth=\z@ \rightwidth=\z@ \topheight=\z@
\botheight=\z@ \setbox\@picbox\hbox\bgroup}
 
\def\enddiagram{\egroup\wd\@picbox\rightwidth\unitlength
\ht\@picbox\topheight\unitlength \dp\@picbox\botheight\unitlength
\hskip\leftwidth\unitlength\box\@picbox}
 
\def\bfig{\begin{diagram}}
\def\efig{\end{diagram}}
\newcount\wideness \newcount\leftwidth \newcount\rightwidth
\newcount\highness \newcount\topheight \newcount\botheight
 
\def\ratchet#1#2{\ifnum#1<#2 \global #1=#2 \fi}
 
\def\putbox(#1,#2)#3{%
\horsize{\wideness}{#3} \divide\wideness by 2
{\advance\wideness by #1 \ratchet{\rightwidth}{\wideness}}
{\advance\wideness by -#1 \ratchet{\leftwidth}{\wideness}}
\vertsize{\highness}{#3} \divide\highness by 2
{\advance\highness by #2 \ratchet{\topheight}{\highness}}
{\advance\highness by -#2 \ratchet{\botheight}{\highness}}
\put(#1,#2){\makebox(0,0){$#3$}}}
 
\def\putlbox(#1,#2)#3{%
\horsize{\wideness}{#3}
{\advance\wideness by #1 \ratchet{\rightwidth}{\wideness}}
{\ratchet{\leftwidth}{-#1}}
\vertsize{\highness}{#3} \divide\highness by 2
{\advance\highness by #2 \ratchet{\topheight}{\highness}}
{\advance\highness by -#2 \ratchet{\botheight}{\highness}}
\put(#1,#2){\makebox(0,0)[l]{$#3$}}}
 
\def\putrbox(#1,#2)#3{%
\horsize{\wideness}{#3}
{\ratchet{\rightwidth}{#1}}
{\advance\wideness by -#1 \ratchet{\leftwidth}{\wideness}}
\vertsize{\highness}{#3} \divide\highness by 2
{\advance\highness by #2 \ratchet{\topheight}{\highness}}
{\advance\highness by -#2 \ratchet{\botheight}{\highness}}
\put(#1,#2){\makebox(0,0)[r]{$#3$}}}

\def\adjust[#1]{} 
 
\newcount \coefa
\newcount \coefb
\newcount \coefc
\newcount\tempcounta
\newcount\tempcountb
\newcount\tempcountc
\newcount\tempcountd
\newcount\xext
\newcount\yext
\newcount\xoff
\newcount\yoff
\newcount\gap%
\newcount\arrowtypea
\newcount\arrowtypeb
\newcount\arrowtypec
\newcount\arrowtyped
\newcount\arrowtypee
\newcount\height
\newcount\width
\newcount\xpos
\newcount\ypos
\newcount\run
\newcount\rise
\newcount\arrowlength
\newcount\halflength
\newcount\arrowtype
\newdimen\tempdimen
\newdimen\xlen
\newdimen\ylen
\newsavebox{\tempboxa}%
\newsavebox{\tempboxb}%
\newsavebox{\tempboxc}%
 
\newdimen\w@dth
 
\def\setw@dth#1#2{\setbox\z@\hbox{\m@th$#1$}\w@dth=\wd\z@
\setbox\@ne\hbox{\m@th$#2$}\ifnum\w@dth<\wd\@ne \w@dth=\wd\@ne \fi
\advance\w@dth by 1.2em}
 
\def\t@^#1_#2{\allowbreak\def\n@one{#1}\def\n@two{#2}\mathrel
{\setw@dth{#1}{#2}
\mathop{\hbox to \w@dth{\rightarrowfill}}\limits
\ifx\n@one\empty\else ^{\box\z@}\fi
\ifx\n@two\empty\else _{\box\@ne}\fi}}
\def\t@@^#1{\@ifnextchar_{\t@^{#1}}{\t@^{#1}_{}}}
\def\to{\@ifnextchar^{\t@@}{\t@@^{}}}
 
\def\t@left^#1_#2{\def\n@one{#1}\def\n@two{#2}\mathrel{\setw@dth{#1}{#2}
\mathop{\hbox to \w@dth{\leftarrowfill}}\limits
\ifx\n@one\empty\else ^{\box\z@}\fi
\ifx\n@two\empty\else _{\box\@ne}\fi}}
\def\t@@left^#1{\@ifnextchar_{\t@left^{#1}}{\t@left^{#1}_{}}}
\def\toleft{\@ifnextchar^{\t@@left}{\t@@left^{}}}
 
\def\two@^#1_#2{\allowbreak
\def\n@one{#1}\def\n@two{#2}\mathrel{\setw@dth{#1}{#2}
\mathop{\vcenter{\lineskip\z@\baselineskip\z@
                 \hbox to \w@dth{\rightarrowfill}%
                 \hbox to \w@dth{\rightarrowfill}}%
       }\limits
\ifx\n@one\empty\else ^{\box\z@}\fi
\ifx\n@two\empty\else _{\box\@ne}\fi}}
\def\tw@@^#1{\@ifnextchar _{\two@^{#1}}{\two@^{#1}_{}}}
\def\two{\@ifnextchar ^{\tw@@}{\tw@@^{}}}
 
\def\tofr@^#1_#2{\def\n@one{#1}\def\n@two{#2}\mathrel{\setw@dth{#1}{#2}
\mathop{\vcenter{\hbox to \w@dth{\rightarrowfill}\kern-1.7ex
                 \hbox to \w@dth{\leftarrowfill}}%
       }\limits
\ifx\n@one\empty\else ^{\box\z@}\fi
\ifx\n@two\empty\else _{\box\@ne}\fi}}
\def\t@fr@^#1{\@ifnextchar_ {\tofr@^{#1}}{\tofr@^{#1}_{}}}
\def\tofro{\@ifnextchar^ {\t@fr@}{\t@fr@^{}}}

\def\mon{\mathop{\m@th\hbox to
      14.6\P@{\lasyb\char'51\hskip-2.1\P@$\arrext$\hss
$\mathord\rightarrow$}}\limits} 
\def\leftmono{\mathrel{\m@th\hbox to
14.6\P@{$\mathord\leftarrow$\hss$\arrext$\hskip-2.1\P@\lasyb\char'50%
}}\limits} 
\mathchardef\arrext="0200       

\def\settypes(#1,#2,#3){\arrowtypea#1 \arrowtypeb#2 \arrowtypec#3}
\def\settoheight#1#2{\setbox\@tempboxa\hbox{#2}#1\ht\@tempboxa\relax}%
\def\settodepth#1#2{\setbox\@tempboxa\hbox{#2}#1\dp\@tempboxa\relax}%
\def\settokens`#1`#2`#3`#4`{%
     \def\tokena{#1}\def\tokenb{#2}\def\tokenc{#3}\def\tokend{#4}}
\def\setsqparms[#1`#2`#3`#4;#5`#6]{%
\arrowtypea #1
\arrowtypeb #2
\arrowtypec #3
\arrowtyped #4
\width #5
\height #6
}
\def\setpos(#1,#2){\xpos=#1 \ypos#2}

\def\settriparms[#1`#2`#3;#4]{\settripairparms[#1`#2`#3`1`1;#4]}%
 
\def\settripairparms[#1`#2`#3`#4`#5;#6]{%
\arrowtypea #1
\arrowtypeb #2
\arrowtypec #3
\arrowtyped #4
\arrowtypee #5
\width #6
\height #6
}
 
\def\resetparms{\settripairparms[1`1`1`1`1;500]\width 500}
 
\resetparms
 
\def\mvector(#1,#2)#3{
\put(0,0){\vector(#1,#2){#3}}%
\put(0,0){\vector(#1,#2){26}}%
}
\def\evector(#1,#2)#3{{
\arrowlength #3
\put(0,0){\vector(#1,#2){\arrowlength}}%
\advance \arrowlength by-30
\put(0,0){\vector(#1,#2){\arrowlength}}%
}}
 
\def\horsize#1#2{%
\settowidth{\tempdimen}{$#2$}%
#1=\tempdimen
\divide #1 by\unitlength
}
 
\def\vertsize#1#2{%
\settoheight{\tempdimen}{$#2$}%
#1=\tempdimen
\settodepth{\tempdimen}{$#2$}%
\advance #1 by\tempdimen
\divide #1 by\unitlength
}
 
\def\putvector(#1,#2)(#3,#4)#5#6{{%
\ifnum3<\arrowtype
\putdashvector(#1,#2)(#3,#4)#5\arrowtype
\else
\ifnum\arrowtype<-3
\putdashvector(#1,#2)(#3,#4)#5\arrowtype
\else
\xpos=#1
\ypos=#2
\run=#3
\rise=#4
\arrowlength=#5
\ifnum \arrowtype<0
    \ifnum \run=0
        \advance \ypos by-\arrowlength
    \else
        \tempcounta \arrowlength
        \multiply \tempcounta by\rise
        \divide \tempcounta by\run
        \ifnum\run>0
            \advance \xpos by\arrowlength
            \advance \ypos by\tempcounta
        \else
            \advance \xpos by-\arrowlength
            \advance \ypos by-\tempcounta
        \fi
    \fi
    \multiply \arrowtype by-1
    \multiply \rise by-1
    \multiply \run by-1
\fi
\ifcase \arrowtype
\or \put(\xpos,\ypos){\vector(\run,\rise){\arrowlength}}%
\or \put(\xpos,\ypos){\mvector(\run,\rise)\arrowlength}%
\or \put(\xpos,\ypos){\evector(\run,\rise){\arrowlength}}%
\fi\fi\fi
}}
 
\def\putsplitvector(#1,#2)#3#4{
\xpos #1
\ypos #2
\arrowtype #4
\halflength #3
\arrowlength #3
\gap 140
\advance \halflength by-\gap
\divide \halflength by2
\ifnum\arrowtype>0
   \ifcase \arrowtype
   \or \put(\xpos,\ypos){\line(0,-1){\halflength}}%
       \advance\ypos by-\halflength
       \advance\ypos by-\gap
       \put(\xpos,\ypos){\vector(0,-1){\halflength}}%
   \or \put(\xpos,\ypos){\line(0,-1)\halflength}%
       \put(\xpos,\ypos){\vector(0,-1)3}%
       \advance\ypos by-\halflength
       \advance\ypos by-\gap
       \put(\xpos,\ypos){\vector(0,-1){\halflength}}%
   \or \put(\xpos,\ypos){\line(0,-1)\halflength}%
       \advance\ypos by-\halflength
       \advance\ypos by-\gap
       \put(\xpos,\ypos){\evector(0,-1){\halflength}}%
   \fi
\else \arrowtype=-\arrowtype
   \ifcase\arrowtype
   \or \advance \ypos by-\arrowlength
       \put(\xpos,\ypos){\line(0,1){\halflength}}%
       \advance\ypos by\halflength
       \advance\ypos by\gap
       \put(\xpos,\ypos){\vector(0,1){\halflength}}%
   \or \advance \ypos by-\arrowlength
       \put(\xpos,\ypos){\line(0,1)\halflength}%
       \put(\xpos,\ypos){\vector(0,1)3}%
       \advance\ypos by\halflength
       \advance\ypos by\gap
       \put(\xpos,\ypos){\vector(0,1){\halflength}}%
   \or \advance \ypos by-\arrowlength
       \put(\xpos,\ypos){\line(0,1)\halflength}%
       \advance\ypos by\halflength
       \advance\ypos by\gap
       \put(\xpos,\ypos){\evector(0,1){\halflength}}%
   \fi
\fi
}
 
\def\putmorphism(#1)(#2,#3)[#4`#5`#6]#7#8#9{{%
\run #2
\rise #3
\ifnum\rise=0
  \puthmorphism(#1)[#4`#5`#6]{#7}{#8}#9%
\else\ifnum\run=0
  \putvmorphism(#1)[#4`#5`#6]{#7}{#8}#9%
\else
\setpos(#1)%
\arrowlength #7
\arrowtype #8
\ifnum\run=0
\else\ifnum\rise=0
\else
\ifnum\run>0
    \coefa=1
\else
   \coefa=-1
\fi
\ifnum\arrowtype>0
   \coefb=0
   \coefc=-1
\else
   \coefb=\coefa
   \coefc=1
   \arrowtype=-\arrowtype
\fi
\width=2
\multiply \width by\run
\divide \width by\rise
\ifnum \width<0  \width=-\width\fi
\advance\width by60
\if l#9 \width=-\width\fi
\putbox(\xpos,\ypos){#4}
{\multiply \coefa by\arrowlength
\advance\xpos by\coefa
\multiply \coefa by\rise
\divide \coefa by\run
\advance \ypos by\coefa
\putbox(\xpos,\ypos){#5} }%
{\multiply \coefa by\arrowlength
\divide \coefa by2
\advance \xpos by\coefa
\advance \xpos by\width
\multiply \coefa by\rise
\divide \coefa by\run
\advance \ypos by\coefa
\if l#9%
   \putrbox(\xpos,\ypos){#6}%
\else\if r#9%
   \putlbox(\xpos,\ypos){#6}%
\fi\fi }%
{\multiply \rise by-\coefc
\multiply \run by-\coefc
\multiply \coefb by\arrowlength
\advance \xpos by\coefb
\multiply \coefb by\rise
\divide \coefb by\run
\advance \ypos by\coefb
\multiply \coefc by70
\advance \ypos by\coefc
\multiply \coefc by\run
\divide \coefc by\rise
\advance \xpos by\coefc
\multiply \coefa by140
\multiply \coefa by\run
\divide \coefa by\rise
\advance \arrowlength by\coefa
\ifcase\arrowtype
\or \put(\xpos,\ypos){\vector(\run,\rise){\arrowlength}}%
\or \put(\xpos,\ypos){\mvector(\run,\rise){\arrowlength}}%
\or \put(\xpos,\ypos){\evector(\run,\rise){\arrowlength}}%
\fi}\fi\fi\fi\fi}}

\newcount\numbdashes \newcount\lengthdash \newcount\increment
 
\def\howmanydashes{
\numbdashes=\arrowlength \lengthdash=40
\divide\numbdashes by \lengthdash
\lengthdash=\arrowlength
\divide\lengthdash by \numbdashes
\increment=\lengthdash
\multiply\lengthdash by 3
\divide\lengthdash by 5
}
 
\def\putdashvector(#1)(#2,#3)#4#5{%
\ifnum#3=0 \putdashhvector(#1){#4}#5
\else
\ifnum#2=0
\putdashvvector(#1){#4}#5\fi\fi}
 
\def\putdashhvector(#1,#2)#3#4{{%
\arrowlength=#3 \howmanydashes
\multiput(#1,#2)(\increment,0){\numbdashes}%
{\vrule height .4pt width \lengthdash\unitlength}
\arrowtype=#4 \xpos=#1
\ifnum\arrowtype<0 \advance\arrowtype by 7 \fi
\ifcase\arrowtype
\or \advance\xpos by 10
    \put(\xpos,#2){\vector(-1,0){\lengthdash}}
    \advance\xpos by 40
    \put(\xpos,#2){\vector(-1,0){\lengthdash}}
\or \advance \xpos by 10
    \put(\xpos,#2){\vector(-1,0){\lengthdash}}
    \advance\xpos by  \arrowlength
    \advance\xpos by  -50
    \put(\xpos,#2){\vector(-1,0){\lengthdash}}
\or \advance\xpos by 10
    \put(\xpos,#2){\vector(-1,0){\lengthdash}}
\or \advance\xpos by \arrowlength
    \advance\xpos by -\lengthdash
    \put(\xpos,#2){\vector(1,0){\lengthdash}}
\or {\advance\xpos by 10
    \put(\xpos,#2){\vector(1,0){\lengthdash}}}
    \advance\xpos by \arrowlength
    \advance\xpos by -\lengthdash
    \put(\xpos,#2){\vector(1,0){\lengthdash}}
\or \advance\xpos by \arrowlength
    \advance\xpos by -\lengthdash
    \put(\xpos,#2){\vector(1,0){\lengthdash}}
    \advance\xpos by -40
    \put(\xpos,#2){\vector(1,0){\lengthdash}}
   \fi
}}
 
\def\putdashvvector(#1,#2)#3#4{{%
\arrowlength=#3 \howmanydashes
\ypos=#2 \advance\ypos by -\arrowlength
\multiput(#1,#2)(0,\increment){\numbdashes}%
    {\vrule width .4pt height \lengthdash\unitlength}
\arrowtype=#4 \ypos=#2
\ifnum\arrowtype<0 \advance\arrowtype by 7 \fi
\ifcase\arrowtype
\or \advance\ypos by \arrowlength \advance\ypos by -40
    \put(#1,\ypos){\vector(0,1){\lengthdash}}
    \advance\ypos by -40
    \put(#1,\ypos){\vector(0,1){\lengthdash}}
\or \advance\ypos by 10
    \put(#1,\ypos){\vector(0,1){\lengthdash}}
    \advance\ypos by \arrowlength \advance\ypos by -40
    \put(#1,\ypos){\vector(0,1){\lengthdash}}
\or \advance\ypos by \arrowlength \advance\ypos by -40
    \put(#1,\ypos){\vector(0,1){\lengthdash}}
\or \advance\ypos by 10
    \put(#1,\ypos){\vector(0,-1){\lengthdash}}
\or \advance\ypos by 10
    \put(#1,\ypos){\vector(0,-1){\lengthdash}}
    \advance\ypos by \arrowlength \advance\ypos by -40
    \put(#1,\ypos){\vector(0,-1){\lengthdash}}
\or \advance\ypos by 10
    \put(#1,\ypos){\vector(0,-1){\lengthdash}}
    \advance\ypos by 40
    \put(#1,\ypos){\vector(0,-1){\lengthdash}}
\fi
}}
 
\def\puthmorphism(#1,#2)[#3`#4`#5]#6#7#8{{%
\xpos #1
\ypos #2
\width #6
\arrowlength #6
\arrowtype=#7
\putbox(\xpos,\ypos){#3\vphantom{#4}}%
{\advance \xpos by\arrowlength
\putbox(\xpos,\ypos){\vphantom{#3}#4}}%
\horsize{\tempcounta}{#3}%
\horsize{\tempcountb}{#4}%
\divide \tempcounta by2
\divide \tempcountb by2
\advance \tempcounta by30
\advance \tempcountb by30
\advance \xpos by\tempcounta
\advance \arrowlength by-\tempcounta
\advance \arrowlength by-\tempcountb
\putvector(\xpos,\ypos)(1,0)\arrowlength\arrowtype
\divide \arrowlength by2
\advance \xpos by\arrowlength
\vertsize{\tempcounta}{#5}%
\divide\tempcounta by2
\advance \tempcounta by20
\if a#8 %
   \advance \ypos by\tempcounta
   \putbox(\xpos,\ypos){#5}%
\else
   \advance \ypos by-\tempcounta
   \putbox(\xpos,\ypos){#5}%
\fi}}
 
\def\putvmorphism(#1,#2)[#3`#4`#5]#6#7#8{{%
\xpos #1
\ypos #2
\arrowlength #6
\arrowtype #7
\settowidth{\xlen}{$#5$}%
\putbox(\xpos,\ypos){#3}%
{\advance \ypos by-\arrowlength
\putbox(\xpos,\ypos){#4}}%
{\advance\arrowlength by-140
\advance \ypos by-70
\ifdim\xlen>0pt
   \if m#8%
      \putsplitvector(\xpos,\ypos)\arrowlength\arrowtype
   \else
   \putvector(\xpos,\ypos)(0,-1)\arrowlength\arrowtype
   \fi
\else
   \putvector(\xpos,\ypos)(0,-1)\arrowlength\arrowtype
\fi}%
\ifdim\xlen>0pt
   \divide \arrowlength by2
   \advance\ypos by-\arrowlength
   \if l#8%
      \advance \xpos by-40
      \putrbox(\xpos,\ypos){#5}%
   \else\if r#8%
      \advance \xpos by40
      \putlbox(\xpos,\ypos){#5}%
   \else
      \putbox(\xpos,\ypos){#5}%
   \fi\fi
\fi
}}
 
\def\putsquarep<#1>(#2)[#3;#4`#5`#6`#7]{{%
\setsqparms[#1]%
\setpos(#2)%
\settokens`#3`%
\puthmorphism(\xpos,\ypos)[\tokenc`\tokend`{#7}]{\width}{\arrowtyped}b%
\advance\ypos by \height
\puthmorphism(\xpos,\ypos)[\tokena`\tokenb`{#4}]{\width}{\arrowtypea}a%
\putvmorphism(\xpos,\ypos)[``{#5}]{\height}{\arrowtypeb}l%
\advance\xpos by \width
\putvmorphism(\xpos,\ypos)[``{#6}]{\height}{\arrowtypec}r%
}}
 
\def\putsquare{\@ifnextchar <{\putsquarep}{\putsquarep%
   <\arrowtypea`\arrowtypeb`\arrowtypec`\arrowtyped;\width`\height>}}
\def\square{\@ifnextchar< {\squarep}{\squarep
   <\arrowtypea`\arrowtypeb`\arrowtypec`\arrowtyped;\width`\height>}}
\def\squarep<#1>[#2`#3`#4`#5;#6`#7`#8`#9]{{
\setsqparms[#1]
\diagram
\putsquarep<\arrowtypea`\arrowtypeb`\arrowtypec`
\arrowtyped;\width`\height>
(0,0)[#2`#3`#4`{#5};#6`#7`#8`{#9}]
\enddiagram
}}                                                 
\def\putptrianglep<#1>(#2,#3)[#4`#5`#6;#7`#8`#9]{{%
\settriparms[#1]%
\xpos=#2 \ypos=#3
\advance\ypos by \height
\puthmorphism(\xpos,\ypos)[#4`#5`{#7}]{\height}{\arrowtypea}a%
\putvmorphism(\xpos,\ypos)[`#6`{#8}]{\height}{\arrowtypeb}l%
\advance\xpos by\height
\putmorphism(\xpos,\ypos)(-1,-1)[``{#9}]{\height}{\arrowtypec}r%
}}
 
\def\putptriangle{\@ifnextchar <{\putptrianglep}{\putptrianglep
   <\arrowtypea`\arrowtypeb`\arrowtypec;\height>}}
\def\ptriangle{\@ifnextchar <{\ptrianglep}{\ptrianglep
   <\arrowtypea`\arrowtypeb`\arrowtypec;\height>}}
\def\ptrianglep<#1>[#2`#3`#4;#5`#6`#7]{{
\settriparms[#1]
\diagram
\putptrianglep<\arrowtypea`\arrowtypeb`
\arrowtypec;\height>
(0,0)[#2`#3`#4;#5`#6`{#7}]
\enddiagram
}}                                            
 
\def\putqtrianglep<#1>(#2,#3)[#4`#5`#6;#7`#8`#9]{{%
\settriparms[#1]%
\xpos=#2 \ypos=#3
\advance\ypos by\height
\puthmorphism(\xpos,\ypos)[#4`#5`{#7}]{\height}{\arrowtypea}a%
\putmorphism(\xpos,\ypos)(1,-1)[``{#8}]{\height}{\arrowtypeb}l%
\advance\xpos by\height
\putvmorphism(\xpos,\ypos)[`#6`{#9}]{\height}{\arrowtypec}r%
}}
 
\def\putqtriangle{\@ifnextchar <{\putqtrianglep}{\putqtrianglep
   <\arrowtypea`\arrowtypeb`\arrowtypec;\height>}}
\def\qtriangle{\@ifnextchar <{\qtrianglep}{\qtrianglep
   <\arrowtypea`\arrowtypeb`\arrowtypec;\height>}}
\def\qtrianglep<#1>[#2`#3`#4;#5`#6`#7]{{
\settriparms[#1]
\width=\height                                
\diagram
\putqtrianglep<\arrowtypea`\arrowtypeb`
\arrowtypec;\height>
(0,0)[#2`#3`#4;#5`#6`{#7}]
\enddiagram
}}
 
\def\putdtrianglep<#1>(#2,#3)[#4`#5`#6;#7`#8`#9]{{%
\settriparms[#1]%
\xpos=#2 \ypos=#3
\puthmorphism(\xpos,\ypos)[#5`#6`{#9}]{\height}{\arrowtypec}b%
\advance\xpos by \height \advance\ypos by\height
\putmorphism(\xpos,\ypos)(-1,-1)[``{#7}]{\height}{\arrowtypea}l%
\putvmorphism(\xpos,\ypos)[#4``{#8}]{\height}{\arrowtypeb}r%
}}
 
\def\putdtriangle{\@ifnextchar <{\putdtrianglep}{\putdtrianglep
   <\arrowtypea`\arrowtypeb`\arrowtypec;\height>}}
\def\dtriangle{\@ifnextchar <{\dtrianglep}{\dtrianglep
   <\arrowtypea`\arrowtypeb`\arrowtypec;\height>}}
\def\dtrianglep<#1>[#2`#3`#4;#5`#6`#7]{{
\settriparms[#1]
\width=\height                                
\diagram
\putdtrianglep<\arrowtypea`\arrowtypeb`
\arrowtypec;\height>
(0,0)[#2`#3`#4;#5`#6`{#7}]
\enddiagram
}}
 
\def\putbtrianglep<#1>(#2,#3)[#4`#5`#6;#7`#8`#9]{{%
\settriparms[#1]%
\xpos=#2 \ypos=#3
\puthmorphism(\xpos,\ypos)[#5`#6`{#9}]{\height}{\arrowtypec}b%
\advance\ypos by\height
\putmorphism(\xpos,\ypos)(1,-1)[``{#8}]{\height}{\arrowtypeb}r%
\putvmorphism(\xpos,\ypos)[#4``{#7}]{\height}{\arrowtypea}l%
}}
 
\def\putbtriangle{\@ifnextchar <{\putbtrianglep}{\putbtrianglep
   <\arrowtypea`\arrowtypeb`\arrowtypec;\height>}}
\def\btriangle{\@ifnextchar <{\btrianglep}{\btrianglep
   <\arrowtypea`\arrowtypeb`\arrowtypec;\height>}}
\def\btrianglep<#1>[#2`#3`#4;#5`#6`#7]{{
\settriparms[#1]
\width=\height                               
\diagram
\putbtrianglep<\arrowtypea`\arrowtypeb`
\arrowtypec;\height>
(0,0)[#2`#3`#4;#5`#6`{#7}]
\enddiagram
}}
 
\def\putAtrianglep<#1>(#2,#3)[#4`#5`#6;#7`#8`#9]{{%
\settriparms[#1]%
\xpos=#2 \ypos=#3
{\multiply \height by2
\puthmorphism(\xpos,\ypos)[#5`#6`{#9}]{\height}{\arrowtypec}b}%
\advance\xpos by\height \advance\ypos by\height
\putmorphism(\xpos,\ypos)(-1,-1)[#4``{#7}]{\height}{\arrowtypea}l%
\putmorphism(\xpos,\ypos)(1,-1)[``{#8}]{\height}{\arrowtypeb}r%
}}
 
\def\putAtriangle{\@ifnextchar <{\putAtrianglep}{\putAtrianglep
   <\arrowtypea`\arrowtypeb`\arrowtypec;\height>}}
\def\Atriangle{\@ifnextchar <{\Atrianglep}{\Atrianglep
   <\arrowtypea`\arrowtypeb`\arrowtypec;\height>}}
\def\Atrianglep<#1>[#2`#3`#4;#5`#6`#7]{{
\settriparms[#1]
\width=\height                                     
\diagram
\putAtrianglep<\arrowtypea`\arrowtypeb`
\arrowtypec;\height>
(0,0)[#2`#3`#4;#5`#6`{#7}]
\enddiagram
}}
 
\def\putAtrianglepairp<#1>(#2)[#3;#4`#5`#6`#7`#8]{{%
\settripairparms[#1]%
\setpos(#2)%
\settokens`#3`%
\puthmorphism(\xpos,\ypos)[\tokenb`\tokenc`{#7}]{\height}{\arrowtyped}b%
\advance\xpos by\height
\puthmorphism(\xpos,\ypos)[\phantom{\tokenc}`\tokend`{#8}]%
{\height}{\arrowtypee}b%
\advance\ypos by\height
\putmorphism(\xpos,\ypos)(-1,-1)[\tokena``{#4}]{\height}{\arrowtypea}l%
\putvmorphism(\xpos,\ypos)[``{#5}]{\height}{\arrowtypeb}m%
\putmorphism(\xpos,\ypos)(1,-1)[``{#6}]{\height}{\arrowtypec}r%
}}
 
\def\putAtrianglepair{\@ifnextchar <{\putAtrianglepairp}{\putAtrianglepairp%
   <\arrowtypea`\arrowtypeb`\arrowtypec`\arrowtyped`\arrowtypee;\height>}}
\def\Atrianglepair{\@ifnextchar <{\Atrianglepairp}{\Atrianglepairp%
   <\arrowtypea`\arrowtypeb`\arrowtypec`\arrowtyped`\arrowtypee;\height>}}
 
\def\Atrianglepairp<#1>[#2;#3`#4`#5`#6`#7]{{
\settripairparms[#1]
\settokens`#2`
\width=\height                                
\diagram
\putAtrianglepairp                            
<\arrowtypea`\arrowtypeb`\arrowtypec`
\arrowtyped`\arrowtypee;\height>
(0,0)[{#2};#3`#4`#5`#6`{#7}]
\enddiagram
}}
 
\def\putVtrianglep<#1>(#2,#3)[#4`#5`#6;#7`#8`#9]{{%
\settriparms[#1]%
\xpos=#2 \ypos=#3
\advance\ypos by\height
{\multiply\height by2
\puthmorphism(\xpos,\ypos)[#4`#5`{#7}]{\height}{\arrowtypea}a}%
\putmorphism(\xpos,\ypos)(1,-1)[`#6`{#8}]{\height}{\arrowtypeb}l%
\advance\xpos by\height
\advance\xpos by\height
\putmorphism(\xpos,\ypos)(-1,-1)[``{#9}]{\height}{\arrowtypec}r%
}}
 
\def\putVtriangle{\@ifnextchar <{\putVtrianglep}{\putVtrianglep
   <\arrowtypea`\arrowtypeb`\arrowtypec;\height>}}
\def\Vtriangle{\@ifnextchar <{\Vtrianglep}{\Vtrianglep
   <\arrowtypea`\arrowtypeb`\arrowtypec;\height>}}
\def\Vtrianglep<#1>[#2`#3`#4;#5`#6`#7]{{
\settriparms[#1]
\width=\height                                 
\diagram
\putVtrianglep<\arrowtypea`\arrowtypeb`
\arrowtypec;\height>
(0,0)[#2`#3`#4;#5`#6`{#7}]
\enddiagram
}}
 
\def\putVtrianglepairp<#1>(#2)[#3;#4`#5`#6`#7`#8]{{
\settripairparms[#1]%
\setpos(#2)%
\settokens`#3`%
\advance\ypos by\height
\putmorphism(\xpos,\ypos)(1,-1)[`\tokend`{#6}]{\height}{\arrowtypec}l%
\puthmorphism(\xpos,\ypos)[\tokena`\tokenb`{#4}]{\height}{\arrowtypea}a%
\advance\xpos by\height
\puthmorphism(\xpos,\ypos)[\phantom{\tokenb}`\tokenc`{#5}]%
{\height}{\arrowtypeb}a%
\putvmorphism(\xpos,\ypos)[``{#7}]{\height}{\arrowtyped}m%
\advance\xpos by\height
\putmorphism(\xpos,\ypos)(-1,-1)[``{#8}]{\height}{\arrowtypee}r%
}}
 
\def\putVtrianglepair{\@ifnextchar <{\putVtrianglepairp}{\putVtrianglepairp%
    <\arrowtypea`\arrowtypeb`\arrowtypec`\arrowtyped`\arrowtypee;\height>}}
\def\Vtrianglepair{\@ifnextchar <{\Vtrianglepairp}{\Vtrianglepairp%
    <\arrowtypea`\arrowtypeb`\arrowtypec`\arrowtyped`\arrowtypee;\height>}}
\def\Vtrianglepairp<#1>[#2;#3`#4`#5`#6`#7]{{
\settripairparms[#1]
\settokens`#2`
\diagram
\putVtrianglepairp                             
<\arrowtypea`\arrowtypeb`\arrowtypec`
\arrowtyped`\arrowtypee;\height>
(0,0)[{#2};#3`#4`#5`#6`{#7}]
\enddiagram
}}

\def\putCtrianglep<#1>(#2,#3)[#4`#5`#6;#7`#8`#9]{{%
\settriparms[#1]%
\xpos=#2 \ypos=#3
\advance\ypos by\height
\putmorphism(\xpos,\ypos)(1,-1)[``{#9}]{\height}{\arrowtypec}l%
\advance\xpos by\height
\advance\ypos by\height
\putmorphism(\xpos,\ypos)(-1,-1)[#4`#5`{#7}]{\height}{\arrowtypea}l%
{\multiply\height by 2
\putvmorphism(\xpos,\ypos)[`#6`{#8}]{\height}{\arrowtypeb}r}%
}}
 
\def\putCtriangle{\@ifnextchar <{\putCtrianglep}{\putCtrianglep
    <\arrowtypea`\arrowtypeb`\arrowtypec;\height>}}
\def\Ctriangle{\@ifnextchar <{\Ctrianglep}{\Ctrianglep
    <\arrowtypea`\arrowtypeb`\arrowtypec;\height>}}
\def\Ctrianglep<#1>[#2`#3`#4;#5`#6`#7]{{
\settriparms[#1]
\width=\height                               
\diagram
\putCtrianglep<\arrowtypea`\arrowtypeb`
\arrowtypec;\height>
(0,0)[#2`#3`#4;#5`#6`{#7}]
\enddiagram
}}                                           
\def\putDtrianglep<#1>(#2,#3)[#4`#5`#6;#7`#8`#9]{{%
\settriparms[#1]%
\xpos=#2 \ypos=#3
\advance\xpos by\height \advance\ypos by\height
\putmorphism(\xpos,\ypos)(-1,-1)[``{#9}]{\height}{\arrowtypec}r%
\advance\xpos by-\height \advance\ypos by\height
\putmorphism(\xpos,\ypos)(1,-1)[`#5`{#8}]{\height}{\arrowtypeb}r%
{\multiply\height by 2
\putvmorphism(\xpos,\ypos)[#4`#6`{#7}]{\height}{\arrowtypea}l}%
}}
 
\def\putDtriangle{\@ifnextchar <{\putDtrianglep}{\putDtrianglep
    <\arrowtypea`\arrowtypeb`\arrowtypec;\height>}}
\def\Dtriangle{\@ifnextchar <{\Dtrianglep}{\Dtrianglep
   <\arrowtypea`\arrowtypeb`\arrowtypec;\height>}}
\def\Dtrianglep<#1>[#2`#3`#4;#5`#6`#7]{{
\settriparms[#1]
\width=\height                              
\diagram
\putDtrianglep<\arrowtypea`\arrowtypeb`
\arrowtypec;\height>
(0,0)[#2`#3`#4;#5`#6`{#7}]
\enddiagram
}}                                          
\def\setrecparms[#1`#2]{\width=#1 \height=#2}%
 
\def\recursep<#1`#2>[#3;#4`#5`#6`#7`#8]{{\m@th
\width=#1 \height=#2
\settokens`#3`
\settowidth{\tempdimen}{$\tokena$}
\ifdim\tempdimen=0pt
  \savebox{\tempboxa}{\hbox{$\tokenb$}}%
  \savebox{\tempboxb}{\hbox{$\tokend$}}%
  \savebox{\tempboxc}{\hbox{$#6$}}%
\else
  \savebox{\tempboxa}{\hbox{$\hbox{$\tokena$}\times\hbox{$\tokenb$}$}}%
  \savebox{\tempboxb}{\hbox{$\hbox{$\tokena$}\times\hbox{$\tokend$}$}}%
  \savebox{\tempboxc}{\hbox{$\hbox{$\tokena$}\times\hbox{$#6$}$}}%
\fi
\ypos=\height
\divide\ypos by 2
\xpos=\ypos
\advance\xpos by \width
\bfig
\putCtrianglep<-1`1`1;\ypos>(0,0)[`\tokenc`;#5`#6`{#7}]%
\puthmorphism(\ypos,0)[\tokend`\usebox{\tempboxb}`{#8}]{\width}{-1}b%
\puthmorphism(\ypos,\height)[\tokenb`\usebox{\tempboxa}`{#4}]{\width}{-1}a%
\advance\ypos by \width
\putvmorphism(\ypos,\height)[``\usebox{\tempboxc}]{\height}1r%
\efig
}}
 
\def\recurse{\@ifnextchar <{\recursep}{\recursep<\width`\height>}}
 
\def\puttwohmorphisms(#1,#2)[#3`#4;#5`#6]#7#8#9{{%
\puthmorphism(#1,#2)[#3`#4`]{#7}0a
\ypos=#2
\advance\ypos by 20
\puthmorphism(#1,\ypos)[\phantom{#3}`\phantom{#4}`#5]{#7}{#8}a
\advance\ypos by -40
\puthmorphism(#1,\ypos)[\phantom{#3}`\phantom{#4}`#6]{#7}{#9}b
}}
 
\def\puttwovmorphisms(#1,#2)[#3`#4;#5`#6]#7#8#9{{%
\putvmorphism(#1,#2)[#3`#4`]{#7}0a
\xpos=#1
\advance\xpos by -20
\putvmorphism(\xpos,#2)[\phantom{#3}`\phantom{#4}`#5]{#7}{#8}l
\advance\xpos by 40
\putvmorphism(\xpos,#2)[\phantom{#3}`\phantom{#4}`#6]{#7}{#9}r
}}
 
\def\puthcoequalizer(#1)[#2`#3`#4;#5`#6`#7]#8#9{{%
\setpos(#1)%
\puttwohmorphisms(\xpos,\ypos)[#2`#3;#5`#6]{#8}11%
\advance\xpos by #8
\puthmorphism(\xpos,\ypos)[\phantom{#3}`#4`#7]{#8}1{#9}
}}
 
\def\putvcoequalizer(#1)[#2`#3`#4;#5`#6`#7]#8#9{{%
\setpos(#1)%
\puttwovmorphisms(\xpos,\ypos)[#2`#3;#5`#6]{#8}11%
\advance\ypos by -#8
\putvmorphism(\xpos,\ypos)[\phantom{#3}`#4`#7]{#8}1{#9}
}}
 
\def\putthreehmorphisms(#1)[#2`#3;#4`#5`#6]#7(#8)#9{{%
\setpos(#1) \settypes(#8)
\if a#9 %
     \vertsize{\tempcounta}{#5}%
     \vertsize{\tempcountb}{#6}%
     \ifnum \tempcounta<\tempcountb \tempcounta=\tempcountb \fi
\else
     \vertsize{\tempcounta}{#4}%
     \vertsize{\tempcountb}{#5}%
     \ifnum \tempcounta<\tempcountb \tempcounta=\tempcountb \fi
\fi
\advance \tempcounta by 60
\puthmorphism(\xpos,\ypos)[#2`#3`#5]{#7}{\arrowtypeb}{#9}
\advance\ypos by \tempcounta
\puthmorphism(\xpos,\ypos)[\phantom{#2}`\phantom{#3}`#4]{#7}{\arrowtypea}{#9}
\advance\ypos by -\tempcounta \advance\ypos by -\tempcounta
\puthmorphism(\xpos,\ypos)[\phantom{#2}`\phantom{#3}`#6]{#7}{\arrowtypec}{#9}
}}
 
\def\setarrowtoks[#1`#2`#3`#4`#5`#6]{%
\def\toka{#1}
\def\tokb{#2}
\def\tokc{#3}
\def\tokd{#4}
\def\toke{#5}
\def\tokf{#6}
}
\def\hex{\@ifnextchar <{\hexp}{\hexp<1000`400>}}
\def\hexp<#1`#2>[#3`#4`#5`#6`#7`#8;#9]{%
\setarrowtoks[#9]
\yext=#2 \advance \yext by #2
\xext=#1 \advance\xext by \yext
\bfig
\putCtriangle<-1`0`1;#2>(0,0)[`#5`;\tokb``\tokd]
\xext=#1 \yext=#2 \advance \yext by #2
\putsquare<1`0`0`1;\xext`\yext>(#2,0)[#3`#4`#7`#8;\toka```\tokf]
\advance \xext by #2
\putDtriangle<0`1`-1;#2>(\xext,0)[`#6`;`\tokc`\toke]
\efig
}
\makeatother

\input BoxedEPS.tex
\SetRokickiEPSFSpecial
\HideDisplacementBoxes

\def\draft{n}
\def\printname#1{
	\if\draft n
		\smash{\makebox[0pt]{\hspace{-0.5in}
			\raisebox{8pt}{\tt\tiny #1}}}
	\fi }

\def\lbl#1{\label{#1}\printname{#1}}
\theoremstyle{plain}

\newtheorem{theorem}{Theorem}
\newtheorem{prop}{Proposition}[section]
\newtheorem{lemma}[prop]{Lemma}
\newtheorem{corollary}{Corollary}
\newtheorem*{claim}{Claim}

\theoremstyle{definition}

\newtheorem{question}{Question}

\theoremstyle{remark}

\newtheorem{remark}[prop]{Remark}

\newcommand{\h}{\operatorname{H_1 }}
\newcommand{\ho}{\operatorname{Hom}}
\newcommand{\im}{\operatorname{Im}}
\newcommand{\aut}{\operatorname{Aut}}
\newcommand{\con}{\equiv}
\newcommand{\rank}{\operatorname{rank}}

\newlength{\standardunitlength}
\def\printname#1{
	\if\draft y
		\smash{\makebox[0pt]{\hspace{-0.5in}
			\raisebox{8pt}{\tt\tiny #1}}}
	\fi
}

\catcode`\@=11
\long\def\@makecaption#1#2{%
    \vskip 10pt
    \setbox\@tempboxa\hbox{
      \small\sf{\bfcaptionfont #1. }\ignorespaces #2}%
    \ifdim \wd\@tempboxa >\captionwidth {%
        \rightskip=\@captionmargin\leftskip=\@captionmargin
        \unhbox\@tempboxa\par}%
      \else
        \hbox to\hsize{\hfil\box\@tempboxa\hfil}%
    \fi}
\font\bfcaptionfont=cmssbx10 scaled \magstephalf
\newdimen\@captionmargin\@captionmargin=2\parindent
\newdimen\captionwidth\captionwidth=\hsize
\catcode`\@=12

\def\L{\Lambda}
\def\Fi{\Phi}
\def\a{\alpha}
\def\b{\beta}
\def\d{\delta}
\def\k{\kappa}
\def\f{\phi}
\def\e{\epsilon}
\def\ss{\sigma}
\def\g{\gamma}
\def\l{\lambda}

\def\fs{f_{\ast}}
\def\fl{F_L}
\def\rl{\hat L}
\def\hj{\hat J}
\def\bj{\bar J}
\def\bh{\bar h}
\def\th{\tilde h}
\def\dg{D_g}
\def\tj{\tilde J}
\def\wj{J_{\omega}}
\def\Gg{{\Gamma_g}}
\def\gg{\Gamma_{g,0}}
\def\LLg{\overline{\mathcal L^L_g}}
\def\llg{\overline{\mathcal L^L_{g,0}}}
\def\hs{homology sphere}
\def\Tg{\mathcal T_g}
\def\tg{\mathcal T_{g,0}}
\def\Kg{\mathcal K_g}
\def\kg{\mathcal K_{g,0}}
\def\Sg{\Sigma_g}
\def\sg{\Sigma_{g,0}}
\def\Lg{\mathcal L^L_g}
\def\lg{\mathcal L^L_{g,0}}
\def\pg{\mathcal P_g}
\def\Z{\mathbb Z}
\def\JJ{\mathcal J}
\def\D{\mathbb L}
\def\F{\mathcal F}
\def\A{\mathcal A}
\def\P{\mathcal P}

\def\S{\mathcal S}
\def\G{\mathcal G}
\def\Q{\mathbb Q}
\def\V{\mathcal V}
\def\M{\mathcal M}
\def\j0{J_0}
\def\s3{S^3}

\def\c{\circ}
\def\iso{\cong}
\def\w{\wedge}
\def\sub{\subseteq}
\def\nsub{\nsubseteq}
\def\sup{\supseteq}
\def\con{\equiv}
\begin{document}
\title{\bf Pure braids, a new
subgroup of the mapping class group and finite-type invariants}
\author{  Jerome Levine  \\\small
        Brandeis University  \\\small Waltham, MA 02254 USA
\\\small levine@binah.cc.brandeis.edu}
\thanks{Partially supported by NSF grant DMS-96-26639}

\date{\normalsize \today }
\maketitle
\pagenumbering{roman}
\begin{abstract} In the study of the relation between the mapping class group
$\Gg$ of a surface $\Sg$ of genus $g$ and the theory of finite-type
invariants of homology $3$-spheres, three subgroups of $\Gg$ play a large
role. They are the Torelli group $\Tg$, the Johnson subgroup $\Kg$ and a new
subgroup $\Lg$, which contains $\Kg$, defined by a choice of a Lagrangian
subgroup $L\sub H_1 (\Sg )$. In this work we determine  the quotient $\Lg
/\Kg\sub\Gg /\Kg$, in terms of the precise description of $\Gg /\Kg$ given by
Johnson and Morita. We also study the lower central series of $\Lg$ and
$\Kg$, using some natural imbeddings of the pure braid group in $\Lg$ and
the theory of finite-type invariants.
\end{abstract}
\tableofcontents
\newpage
\pagenumbering{arabic}
\section{Introduction}

Let
$\Gg$ denote the group of isotopy classes of  orientation-preserving
diffeomorphisms of the {\em bounded }surface
$\Sg$ of genus
$g$ with one boundary component which are the identity on the boundary.
In~\cite{GL} and~\cite{GL1} the relation between the mapping class group
and the
theory of finite-type invariants of
\hs s is investigated. It turns out that the three different notions of
finite-type which are discussed in~\cite{GL} require one to consider three
different subgroups of the mapping class group. Two of these subgroups are
familiar from previous work of D. Johnson,
S. Morita and others on the structure of the mapping class group
(see~\cite{J2}) and~\cite{M1} for surveys). These are the Torelli group
$\Tg$ (maps which induce the identity automorphism on $\h (\Sg )$)  and the
{\em Johnson subgroup } $\Kg$ (generated by Dehn twists on bounding closed
curves). But the original notion of finite-type introduced by Ohtsuki
in~\cite{O} requires one to consider a new subgroup defined as follows.
Choose a
{\em Lagrangian } subspace $L\sub H=\h  (\Sg )$, i.e. a summand of rank $g$
on which the intersection form vanishes. Define $\Lg$ to be the subgroup of
$\Gg$ generated by Dehn twists on closed curves whose homology class lies in
$L$. Of course $\Lg$ depends on $L$ but its conjugacy class is independent
of the choice of $L$. It is also clear that $\Kg\sub\Lg$. In fact, the
relation between $\Lg$ and $\LLg$ is comparable to the relation between $\Kg$
and
$\Tg$.

In this work we will give a precise determination of how $\Lg$ sits in
$\Gg$, or, more precisely, how $\Lg /\Kg$ sits in $\Gg /\Kg$, using the
description given by Johnson~\cite{J} and~\cite{J1}. The analogous result
for a {\em closed }surface is obtained as a corollary and is, in fact,
somewhat simpler to state than the bounded case. 

A second purpose of this paper is to study the lower central series of the
groups $\Lg$ and $\Kg$. We show that the lower central series of $\Lg$ is
dominated by a refinement of the relative weight filtration of $\Gg$. We also
make use of  a natural imbedding of the pure braid group $P_g$ on $g$
strands into $\Lg$ (first used by Oda, in an unpublished paper
\cite{Od}\footnote{The author thanks S. Morita for informing him of the
existence of this paper and providing a copy.} to give
lower bounds for the ranks of the associated gradeds of the relative weight
filtration of $\Gg$), and of $P_{2g}$ into $\Kg$, to give lower bounds on
the ranks of the lower central series quotients of these groups. These estimates
augment the lower bounds given in
\cite{GL1} using the theory of finite-type invariants of homology $3$-spheres.

 Of course the definitive results on the lower
central series of $\Tg$ are those of Hain \cite{Ha}.

\section{Statement of Results}
\subsection{The group $\Lg$}\lbl{sec.L} Let $\LLg$ denote the subgroup of
$\Gg$ consisting of maps whose induced automorphism of $H$ is the identity
on $L$. Then $\Lg\sub\LLg$ and
$\Tg\sub\LLg$.  Under the obvious
isomorphism
$\Gg /\Tg
\iso Sp(\Z ,2g)$, the symplectic group over $\Z$ of genus $g$, $\LLg /\Tg$
corresponds to the subgroup
$B_g$ consisting of symplectic automorphisms of $H$ which are the identity on
$L$. $B_g$ is isomorphic to $\S_2 (L)$, the additive subgroup of symmetric
elements of $L\otimes L$. If $f\in B_g$, then
$\f-1$ induces an element of $\ho(H/L,L)\iso\ho(L^{\ast},L)\iso  L\otimes
L$ (using duality via the intersection pairing on $H$). Alternatively if we
choose a symplectic basis
$\{ x_i ,y_i\}$ for $H$, where $\{x_i\}$ is a basis of $L$, then $\f$ is
represented by a matrix in the form $(\smallmatrix I & A \\
 0 & I \endsmallmatrix )$
 and
$\f\to A$ defines an isomorphism of $B_g$ with the additive group of
symmetric integral matrices. See~\cite{GL} for more details.

Recall the Johnson homomorphism $J:\Tg\to\L^3 H$, which is onto and whose
kernel is the Johnson subgroup $\Kg$ (see~\cite{J} and~\cite{J1}). If
$f\in\Tg$ and
$\a\in F=\pi_1 (\Sg )$, the free group on $2g$ generators, then we can write
\begin{displaymath} f_{\ast}(\a )=\a \Fi_f (\a ) \bmod F_3
\end{displaymath} where, for any group $G$, $G_q$ denotes the $q$-th lower
central series subgroup: $G_1 =G,\ G_{q+1}=[G,G_q ]$. $\Fi_f (\a )$ can be
regarded, therefore, as an element in $F_2 /F_3 \iso\L^2 H$ and so
$\Fi_f$ defines a homomorphism $H\to\L^2 H$ or, dually, an element of
$H\otimes\L^2 H$. Johnson shows in~\cite{J} that $\Fi_f$ lies in the subgroup
$\L^3 H\sub H\otimes \L^2 H$, where this imbedding is defined by
\begin{displaymath} h_1 \w h_2\w h_3\to h_1 \otimes (h_2\w h_3 ) + h_2\otimes
(h_3\w h_1  ) + h_3\otimes (h_1 \w h_2 )
\end{displaymath} and, in~\cite{J1}, that the correspondence $f\to\Fi_f$
defines an {\em isomorphism } $\Tg /\Kg\iso
\L^3 H$. From this we can, therefore, regard $\LLg  /\Kg$ as an extension of
$\L^3 H$ by $B_g$.
\medskip
\begin{theorem}\lbl{th.L} There exists a homomorphism
\begin{equation}\lbl{eq.J}
\JJ :\LLg\to\L^3 (H/L)\oplus H/L
\end{equation} extending the Johnson homomorphism, which satisfies the
following:
\begin{enumerate}
\item[(a)] $\JJ$ is onto,
\item[(b)] $\ker\JJ =\Lg$
\item[(c)] There is a commutative diagram
\begin{equation}\lbl{eq.John}
\bfig
\putsquare<1`-1`-1`1;500`500>(500,500)[\Lg`\phantom{\LLg}`\Kg`\phantom{\Tg};```]
\putsquare<1`-1`-1`1;700`500>(1000,500)[\LLg `\L^3 (H/L)\oplus H/L `\Tg
`\L^3 H ;\JJ
` `p\oplus c `J]
\putmorphism(1000,1500)(0,-1)[B_g`\phantom{\LLg}`]{500}{-1}m
\putmorphism(0,1000)(1,0)[1`\phantom{\Lg}`]{500}1a
\putmorphism(0,500)(1,0)[1`\phantom{\Kg}`]{500}1a
\putmorphism(500,500)(0,-1)[\phantom{\Kg}`1`]{500}{-1}m
\putmorphism(1000,500)(0,-1)[\phantom{\Tg}`1`]{500}{-1}m
\putmorphism(1700,1500)(0,-1)[1`\phantom{\L^3 (H/L)\oplus H/L}`]{500}{-1}m
\putmorphism(1700,1000)(1,0)[\phantom{\L^3 (H/L)\oplus H/L}`1`]{600}{1}a
\putmorphism(1700,500)(1,0)[\phantom{\L^3 H}`1`]{600}{1}a
\efig
\end{equation}
 where $p$ is the projection $\L^3 H\to\L^3 (H/L)$ and
$c$ is the composition of the projection $H\to H/L$ with the contraction map
$\L^3 H\to H$ defined by
$$ h_1\w h_2\w h_3\to (h_1\cdot h_2 )h_3 + (h_2\cdot h_3 )h_1  + (h_3\cdot
h_1  )h_2
$$
\end{enumerate}
\end{theorem}
\medskip

If $\gg$ denotes the mapping class group of a {\em closed }surface
$\sg$ of genus $g$, i.e. the group of isotopy clases of
orientation-preserving diffeomorphisms of $\sg$, then we can define the
analogous subgroups $\llg, \lg\sub\gg$. We also have the corresponding
Torelli group $\tg$ and Johnson subgroup $\kg$. The homomorphism
$J:\Tg\to\L^3 H$ is shown, in~\cite{J}, to induce a homomophism
$J':\tg\to\L^3 H/H$ whose kernel is $\kg$.

\begin{corollary}\lbl{cor.cl} There exists an epimorphism $\j0 :\llg\to\L^3
(H/L)$ whose kernel is $\lg$ and which satisfies the commutative diagram:
\begin{equation}
\begin{CD}
1 @>>>  \lg @>>> \llg @>\j0 >> \L^3 (H/L) @>>> 1 \\
@.  @AAA @AAA @A p' AA @. \\
1 @>>> \kg @>>> \tg @> J' >> \L^3 H/H @>>> 1
\end{CD}
\end{equation}
where $p'$ is induced by $p$.
\end{corollary}
We can also characterize $\lg\sub\gg$ in terms of the induced outer
automorphisms
of $\pi =\pi_1 (\sg )$. let $\pi_L\sub\pi$ be the subgroup of elements whose
homology class belongs to $L$. If $f\in\llg$, then
$f_{\ast}|\pi_L\con\text{id}\mod\pi_2$.
\begin{corollary}\lbl{cor.pl}
If $f\in\gg$, then $f\in\lg$ if and only if $f_{\ast}|\pi_L\con\text{id}\mod
[\pi ,\pi_L ]$.
\end{corollary}
Note that this condition on $f_{\ast}$ depends only on its outer
automorphism class.

 We also identify the lower central
series terms of
$\Lg$, modulo
$\Kg$. Let
\mbox{$K=\ker (p\oplus e)$}. By Theorem~\ref{th.L}, we have an isomorphism
\begin{equation}\lbl{eq.J1}
\JJ':\Lg\cap\Tg /\Kg\iso K
\end{equation}
 We now define a filtration of $\L^3 H$:
$$ \L^3 H=K_0\sup K_1\sup K_2\sup K_3 =\L^3 L\sup K_4 =0  $$ by letting
$K_m$ be the subgroup of $\L^3 H$ generated by elements
$h_1 \w h_2\w h_3$, where at least $m$ of the $h_i$ lie in $L$. Note that
$K_1\supsetneq K\supsetneq K_2$.
\medskip
\begin{theorem}\lbl{th.lcL} Under the isomorphism $\JJ'$ of~\eqref{eq.J1}
the lower central series terms of $\Lg$ satisfy:
$$ ((\Lg )_m\cdot\Kg )/\Kg \iso K_m \text{  for }m\ge 2   $$ In
particular
$(\Lg )_4\sub\Kg$, but $(\Lg )_3\nsubseteq\Kg$.\newline
\vspace{.01mm}

Furthermore $\Lg\cap\Tg =\Kg\cdot [\Lg ,\Tg ]$.
\end{theorem}
\medskip
\begin{question}  Does $[\Lg ,\Lg ]$ contain $\Kg$?
\end{question}

We obtain some further information on the lower central series of $\Lg$
inside
$\Kg$. Recall the filtration on $\Gg$ defined by:
$$ \Gg [n]=\{ f\in\Gg | f_{\ast}:F\to F \text{ is the identity mod }
F_{n+1}\}
$$
This is called the {\em relative weight filtration }in~\cite{Od}. Thus $\Gg
[1]=\Tg$ and
$\Gg [2]=\Kg$.
\medskip
\begin{theorem}\lbl{th.lcl}
$(\Lg )_q \sub\Gg [n]$ if $q\ge\binom {n+3}2 -6$
\end{theorem}  Thus $\Lg$ is residually nilpotent.
\medskip

In fact we have a more delicate description of the relation between the lower
central series of $\Lg$ and the relative weight filtration which requires us
to define a more refined filtration. We first
recall the extended Johnson homomorphisms from~\cite{J2} as formulated by
Morita~\cite{M}. Let
$\D (H)$ be the graded free Lie algebra (over $\Z$) on $H$. $\D_n
(H)$ is generated by the brackets of length
$n$. Then
$$ J_n :\Gg [n]\to H\otimes\D_{n+1}(H) $$ is defined by $f_{\ast}(\a )=\a
J_n (\a ) \bmod F_{n+2}$, where $J_n (\a )\in F_{n+1}/F_{n+2}
\iso\D_{n+1}(H)$. Then
$\Gg [n+1]=\ker J_n$ and
\begin{equation}\lbl{eq.b}
\im
J_n\sub\ker\{b:H\otimes\D_{n+1}(H)\to\D_{n+2}(H)\}
\end{equation}
 where $b(h\otimes\a
)=[h,\a ]$ (see~\cite{M3}). Thus $J_1 ((\Lg )_m )=K_m$ (by Theorem
\ref{th.lcL}).

We generalize the filtration $\{ K_m\}$ as follows:
$$ \D_m (H)=\D_m^0 (H)\sup\D_m^1 (H)\sup\cdots\sup\D_m^m (H)=\D_m
(L)\sup\D_m^{m+1}(H)=0  $$ where $\D_m^r (H)$ is the subgroup generated by
brackets of length $m$ of which at least $r$ of the entries belong to
$L$. We can then define a filtration of $H\otimes\D_m (H)$ by $\F_m^r
=(L\otimes\D_m^{r-1}(H))\oplus (H\otimes\D_m^r (H))$. Thus
$\F_m^{m+1}=L\otimes\D_m (L)$ and $\F_m^{m+2}=0$.
\medskip
\begin{theorem}\lbl{th.lcsL} $J_n ((\Lg )_{q+r})\sub\F_{n+1}^r$ if
$q\ge\binom{n+3}2 -6$.
\end{theorem}
\medskip

\subsection{Imbedding the pure braid group in the mapping class group} It
has been
pointed out several times in the literature that one can define interesting
maps from
the group of {\em framed} pure braids to the mapping class group. The
framed pure
braid group on $h$ strands, $\P_h$, which is canonically isomorphic to the
Cartesian
product of the usual pure braid group $\P_h^{\c}$ with the free abelian
group of
rank $h$, can be represented as the group of isotopy classes of
diffeomorphisms of
the disk with $h$ holes, $D_h$, which are the identity on the boundary. Any
imbedding
of
$D_h$ into
$\Sg$ defines, therefore, a homomorphism
$\theta:\P_h\to\Gg$. For example one can imbed $D_{2g}$ into $\Sg$ as in Figure
\ref{fig4}.
\begin{figure}[ht]
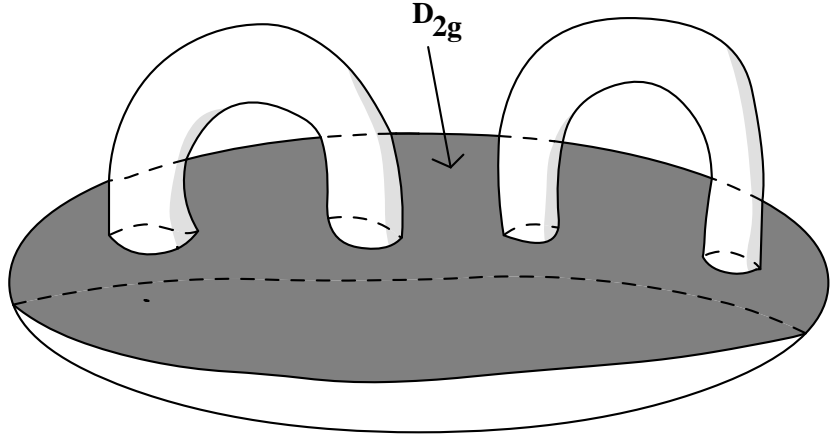

 \centerline{\BoxedEPSF{Fig4 scaled 800}}
\bigskip
\caption{The Hatcher-Thurston map}\lbl{fig4}
\end{figure}
This map was considered by Hatcher-Thurston \cite{HT} and, with a
formulation in
terms of tangles, Matveev-Polyak \cite{MP}  in a slightly more general
setting of
what Matveev-Polyak call {\em admissible} braids. It is not hard to see
that $\theta
(\P_{2g})\sub\Lg$, where
$L$ is the subspace spanned by the meridians of the handles.
We will be more interested, though, in a smaller version of this map.
Consider the
imbedding $D_g\sub\Sg$ indicated in Figure \ref{fig1}.
\begin{figure}[ht]
 \centerline{\BoxedEPSF{Fig1 scaled 800}}
\bigskip
\caption{Imbedding of $\pg$ into $\Lg$}\lbl{fig1}
\end{figure}
The resulting map $\psi
:\pg\to\Gg$ was considered in~\cite{Od} Oda who showed that $\psi$ induces {\em
imbeddings }of the lower-central series quotients of $\pg$ into the
associated graded
quotients of
$\Gg$ using  the relative weight filtration: $\psi_n :(\pg )_n /(\pg
)_{n+1}\to\Gg [n-1]/\Gg [n]$. The well-known split exact sequence $1\to
F^{g-1}\to\P_g\to\P_{g-1}\to 1$ gives an isomorphism:
\begin{equation}\lbl{eq.pblc}
(\pg )_n /(\pg )_{n+1}\iso (\P_{g-1})_n /(\P_{g-1})_{n+1}\oplus F^{g-1}_n
/F^{g-1}_{n+1}
\end{equation}
See e.g.~\cite{F}.This enables Oda to use $\psi_n$ to give some explicit
lower bounds for the rank of $\Gg [n-1]/\Gg [n]$.

We observe easily that $\psi$ actually maps into $\Lg$, for the correct
choice of
$L$, and describe how the image fits into the refined filtration $\F_m^r$ of
$H\otimes\D (H)$.
\medskip
\begin{theorem}\lbl{th.br}
$\psi (\pg )\sub\Lg$ and $\im (J_{n-1}\psi_n )\sub\F_n^{n+1}=L\otimes\D_n
(L)$.
\end{theorem}
\medskip

Thus $(\pg )_n /(\pg )_{n+1}\sub\Gg [n-1]/\Gg [n]$ lies in the
{\em bottom stage } of this filtration.
\begin{remark}\lbl{rem.1}
It is not hard to prove, by a similar argument, that the conclusions of Theorem
\ref{th.br} are also true for the Hatcher-Thurston map.
\end{remark}

As in Equation~\eqref{eq.b} $\im (J_{n-1}\psi_n )\sub\ker\{ b:L\otimes\D_n
(L)\to\D_{n+1}(L)\}$, but for $n\ge 3$, this inclusion is definitely
proper. One sees this by computing the difference between the rank of
$\ker b$ and the rank of $(\pg )_n /(\pg )_{n+1}$, using
Equation~\eqref{eq.pblc}.
For example, for
$n=3$ this difference is $\frac 16 (g^3 -g)$ and for $n=4$ it is $\frac 18
(g^3 -g)(g-2)$. For $n\le 2$ see Theorem~\ref{th.br1} below.

The non-triviality of $J_{n-1}\psi_n$  shows that $\Lg$ is not
nilpotent. More precisely, combining this with Theorem~\ref{th.lcl} gives a
lower bound on the rate of descent of the lower central
series of $\Lg$. We denote the rank of $(\pg )_n /(\pg )_{n+1}$ by $r(g,n)$. 
This can be computed explicitly from Equation~\eqref{eq.pblc}.
\medskip
\begin{corollary}\lbl{cor.p} The image of the map
$(\pg )_n
/(\pg )_q\to (\Lg )_n /(\Lg )_q$ has rank $\ge r(g,n)$ if
$q\ge\binom {n+3}2 -6$.
\end{corollary}
 This estimate is independent of the estimate on the rank of the lower central
series of $\Lg$ given in \cite{GL1}. Recall the map $\phi_n^L :(\Lg )_{3n}/(\Lg
)_{3n+1}\otimes\Q\to\A_n^{\text{conn}}(\emptyset )$ defined in \cite{GL1} and
shown there (Theorem 7) to be onto if $g\ge 5n+1$, where
$\A_n^{\text{conn}}(\emptyset )$ is a vector space defined by trivalent graphs
with $2n$ vertices and $3n$ edges. Now it is easy to see that the composition
$(\pg )_{3n}\to  (\Lg )_{3n}\to\A_n^{\text{conn}}(\emptyset )$ is zero. If
$h\in\im\psi$ then $S^3_h =S^3$ since $h$ extends to a diffeomorphism of the
handlebody bounded by $\Sg$. Thus combining Corollary~\ref{cor.p} and
\cite[Theorem
7]{GL1} we have:
\begin{corollary} $\rank (\Lg )_{3n} /(\Lg )_q\ge
r(g,3n)+\dim\A_n^{\text{conn}}(\emptyset )$ if $q\ge\binom {3n+3}2 -6$ and
$g\ge
5n+1$.
\end{corollary}
\medskip

For the special cases $n=1,2$ we have:
\medskip
\begin{theorem}\lbl{th.br1}
 \begin{enumerate}
\item[(a)] $\psi_1$ induces an isomorphism $\pg /(\pg )_2 \iso\Lg
/(\Lg\cap\Tg )\iso B_g.$
\item[(b)] $J_1\psi_2$ induces an isomorphism $(\pg )_2 /(\pg )_3 \otimes\Q
\iso ((\Lg )_3\cdot\Kg )/\Kg \otimes\Q\iso\L^3 L\otimes\Q$. Thus
$\psi^{-1}(\Kg )=(\pg )_3$.
\end{enumerate}
\end{theorem}
\medskip
\begin{question} Can (b) be improved to say $(\pg )_2 /(\pg )_3\iso\L^3 L$?
\end{question}
\begin{remark}
For the Hatcher-Thurston map $\theta$ one has $J_1\theta ((\P_{2g})_2 )\sub K_3
=\L^3 L$, according to Remark \ref{rem.1}, but it is also not hard to show
that $J_1
(\theta (\P_{2g})\cap\Tg )\sub K_2$. A consequence of this is that $\theta
(\P_{2g})\subsetneq\Lg$.
\end{remark}

 We will consider one more mapping of the pure braid group
into the mapping class group. Consider the imbedding
$D_g\sub\Sg$ indicated in Figure
\ref{fig5}.
\begin{figure}[ht]
 \centerline{\BoxedEPSF{Fig5 scaled 800}}
\bigskip
\caption{Imbedding of $\pg$ into $\Kg$}\lbl{fig5}
\end{figure}
Denote
the resulting map by \mbox{$\k :\pg\to\Gg$.} We will prove:

\begin{theorem}\lbl{th.k}
$\k (\pg)\sub\Kg$ and $\k ((\pg )_n )\sub\Gg [2n]$. The induced map $\k_n :(\pg
)_n /(\pg )_{n+1}\to\Gg [2n]/\Gg [2n+1]$ is a monomorphism. Thus the
induced map $(\pg
)_n /(\pg )_{n+1}\to (\Kg )_n /(\Kg )_{n+1}$ is also a monomorphism and so
$\rank (\Kg )_n /(\Kg )_{n+1}\ge r(g,n)$.
\end{theorem}

As above we have, from \cite{GL1}, a map $\phi_n^K :(\Kg )_n /(\Kg
)_{n+1}\to\A_n^{\text{conn}}(\emptyset )$ which is onto if $g\ge 5n+1$.
The composition $(\pg
)_n /(\pg )_{n+1}\to (\Kg )_n /(\Kg )_{n+1}\to\A_n^{\text{conn}}(\emptyset )$
is zero since, if $h\in\im\phi_n^K$ then $S^3_h =S^3$. Thus we have, from
Theorem~\ref{th.k} and \cite{GL1}:
\begin{corollary}  $\rank(\Kg )_n /(\Kg )_{n+1}\ge r(g,n)+\dim\A_n^{\text{conn}}(\emptyset )$
 if $g\ge 5n+1$.
\end{corollary}

\subsection{Applications of the theory of finite-type invariants}
The fact that there are useful connections between the theory of finite-type
invariants of \hs s and the structure of subgroups of the mapping class
group was
established in \cite{GL} and \cite{GL1}. In particular, as we have mentioned
above, Theorem 6 of
\cite{GL1} gives
lower bounds for the dimensions of the graded quotients of the lower
central series
of the subgroups $\Lg , \Tg$ and $\Kg$ in terms of certain vector spaces
$\A^{\text{conn}}_n (\emptyset )$ defined by connected trivalent graphs.

As a further application of this result we have the following
relation between the lower central series of these groups.
\medskip
\begin{theorem}\lbl{th.lkt}
 Let $g\ge 5n+1$. Then
\begin{enumerate}
\item $(\Kg )_n\nsub (\Tg )_{2n+1}\cup
(\Lg )_{3n+1}$,
\item $(\Tg )_{2n}\nsub (\Kg )_{n+1}\cup (\Lg
)_{3n+1}$,
\item $(\Lg )_{3n}\nsub (\Tg )_{2n+1}\cup (\Kg
)_{n+1}$.
\end{enumerate}
\end{theorem}
\medskip

We also examine the map
$\psi :\pg\to\Lg$, defined in the previous section, and its relationship to the
theory of finite-type invariants of knots and
\hs s. Let
$\V$ denote the
$\Q$-vector space generated by isotopy classes of oriented smooth knots in
$\s3$ and let $\M$ denote the $\Q$-vector space generated by diffeomorphism
classes of smooth oriented \hs s. In~\cite{G} Garoufalidis defines a natural
linear map
$\Psi :\V\to\M$  by $\Psi (K)=\s3_K$, where, for any knot $K$ in a
\hs\ $M$, $M_K$ denotes the \hs\  obtained by doing a $+1$-surgery   on
$K\sub M$. Filtrations are defined on $\V$ and $\M$ by the following
devices. For any positive integer
$q$, a {\em singular
$q$-knot } is an immersion of $S^1$ into $\s3$ with exactly $q$ transverse
(ordered) double points. If $K$ is a singular $q$-knot, the resolution
$\hat K$ of $K$ is the element of $\V$ defined by $\hat K =\sum_{\e}(-1)^{|\e
|}K_{\e}$, where $\e$ ranges over all sequences
$(\e_1 ,\ldots ,\e_q ), \e_i =0$ or $1,\ |\e |=\sum_i \e_i$ and $K_{\e }$
denotes the knot obtained from $K$ by resolving each double point to give a
positive, resp. negative crossing according to whether $\e_i$ is, resp.
$0$ or
$1$. Now define $\F_q (\V )$ to be the subspace of $\V$ generated by the
resolutions of all $q$-singular knots. This gives a decreasing filtration of
$\V$ (see~\cite{BN}). To deal with $\M$ define, for any
$q$-component link $J$ in a \hs , an element $[M;J]\in\M$ by the formula
$[M;J]=\sum_{\e}(-1)^{|\e |}M_{J_{\e}}$ where $\e$ is as above and
$J_{\e}$ is the sublink of $J$ consisting of exactly those components for
which
$\e_i =1$. Now $\F_q (\M )$ is defined to be the subspace of
$\M$ generated by all $[M;J]$ for which $J$ has $q$ components. This is a
decreasing filtration of $\M$ (see~\cite{O}). It is shown in~\cite{O}
and~\cite{GL2} that
$\F_{3q-2}(\M )=\F_{3q}(\M )$, for all $q$.

The associated graded quotients
$$\G_q (\V )=\F_q (\V )/\F_{q+1}(\V )\text{  and }
\G_q (\M )=\F_{3q}(\M )/\F_{3q+1}(\M )$$
 are described in terms of spaces of trivalent graphs $\A (S^1 ), \A
(\emptyset )$, where $\A (X)$, for any
$1$-manifold $X$, is a certain $\Q$-vector space defined from trivalent
graphs which contain $X$ as a prescribed subgraph (see~\cite{BN},~\cite{O}
and~\cite{LMO}). It was conjectured by Garoufalidis in~\cite{G} and proved by
Habegger~\cite{H} that
$\Psi (\F_{2q}(\V ))\sub\F_{3q}(\M )$ and so we have induced maps $\Psi_q
:\G_{2q}(\V )\to\G_q (\M )$. This map is studied in~\cite{GH}, using the work
 of~\cite{LMO}, and shown to be non-trivial for all
$q$.

The map $\psi :\pg\to\Lg$ and $\Psi :\V\to\M$ are related by means of
\lq actions\rq\  of $\pg^{\c}$ on $\V$ and $\Lg$ on $\M$. Let $K$ be an
oriented
knot in $\s3$ and $D\sub\s3$ an \lbl{action} imbedded disk which meets
$K$ transversely in $g$ points so that, at each intersection point,
$K$ points in the same normal direction to $D$. If $b$ is any
$g$-strand pure braid then we can cut $\s3$ along $D$ and insert
$b$ so that the orientations agree, to form a new knot $K_b$. This action was
defined by Stanford in~\cite{S} who showed that if $b\in (\pg^{\c})_m$, then
$K_b -K\in\F_m (\V )$.  Note also that if $b$ is a framed pure braid and $K$ is
a framed knot, then $K_b$ is also framed in a natural way. If $s$ is the
self-linking number of $K$ and $\{ s_i\}$ are the self-linking numbers of the
strands of $b$, then the self-linking number of $K_b$ is $s+\sum_i s_i$.

We now
describe the action of $\Lg$ on
$\M$. Suppose
$M$ is a \hs\ and $i:\Sg\sub M$ is an imbedding. If $L_i =L\sub H=\h (\Sg
)$ is the
kernel of the inclusion into one of the complementary summands, then, for any
$f\in\Lg$, we can cut
$M$ open along
$i(\Sg )$ and reglue, using $f$, to define a new homology sphere $M_f$
(see~\cite{GL}). It is shown in~\cite{GL} that if $f\in (\Lg )_m$ then $M_f
-M\in\F_m (\M )$.

Now suppose we are given $K,D$ as above and $b\in\pg^{\c}$. Choose an imbedding
$i:\Sg\to\s3$ whose image is the boundary of a regular neighborhood of
$K\cup D$. Then $i$ will also define an imbedding $j$ of $\Sg$ into
$\s3_K$. Let $b$ denote also the framed braid whose self-linking numbers are
determined by the equations $\sum_j l_{ij}=0$, where $l_{ij}$ denotes the
linking
number of the $i$-th and $j$-th strand, if $i\not= j$, and the self-linking
number of the $i$-th strand if $i=j$. Then for every
$i=1,\cdots g$, then we have:
\medskip
\begin{theorem}\lbl{th.ft}
$\Psi (K_b )=M_{\psi (b)}$, where $M=\s3_K$.
\end{theorem}
\medskip

\medskip

\section{Proof of Theorem~\ref{th.L}} As mentioned already, different
choices for
$L$ result in conjugate subgroups
$\Lg$ (and also $\LLg$). Therefore it will suffice to prove
Theorems~\ref{th.L}-\ref{th.lcsL} for any convenient choice of $L$. (The
other theorems will require certain choices for $L$). In particular we will
choose
$L$ to be the subgroup of $H$ generated by meridians of the handles for some
representation of $\Sg$ as the boundary of a handle-body. Thus we can choose
a basis $\{ x_i ,y_i\}$ for $F=\pi_1 (\Sg )$ to satisfy:
\begin{itemize}
\item The induced basis of $H$ (which we will also denote by $x_i ,y_i$) is
{\em symplectic },i.e. under the intersection pairing on $H$, we have
$x_i\cdot x_j =y_i\cdot y_j =0$ and $x_i\cdot y_j =\delta_{ij}$,
\item $\{ x_i\}$ is a basis for $L$,
\item The element $[x_1 ,y_1 ]\cdots [x_g ,y_g ]\in F$ represents the
boundary curve of $\Sg$.
\end{itemize} We will call such a basis {\em admissible for $L$}, We will
also use this term for the induced basis of $H$.
\subsection{Construction of $\JJ$}\lbl{sec.J} We begin by defining a crossed
homomorphism on $\LLg$. Let $\fl\sub F$ be the set of all elements of $F$
which map into $L\sub H$ under abelianization. Let
$\rl$ be the abelianization of $\fl$. We define
$$ \hj :\LLg\to\ho(\rl ,\L^2 H) $$ by the formula $\hj
(f)\b\con\a^{-1}f_{\ast}(\a )\ \bmod F_3$, under the identification $\F_2
/F_3\iso\L^2 H$, where $\a\in\fl$ is any lift of
$\b\in\rl$.

We then have the following crossed-homomorphism formula
\begin{equation}\lbl{eq.cross}
\hj (fg)=\hj (g)+\hj (f)\circ g_{\ast}
\end{equation} where $g_{\ast}$ is the automorphism of $\rl$ defined by
$g$.

Composing $\hj$ with the projection $\L^2 H\to\L^2 (H/L)$ defines
$$\bj :\LLg\to\ho(\rl ,\L^2 (H/L)) $$ Composing $\hj$ with the contraction
map
$c:\L^2 H\to\Z$ defines
$$\wj :\LLg\to\ho(\rl ,\Z ) $$
$c$ is defined, using the intersection pairing on $H$, by $h_1 \wedge h_2\to
h_1 \cdot h_2$.

 We point out two properties of $\bj$ and $\wj$.
\begin{lemma}\lbl{lem.cross}
\begin{enumerate}
\item[(a)] $\bj (f)|F_2 =\wj (f)|F_2 =0$, for any $f\in\LLg$,
\item[(b)] $\bj$ and $\wj$ are homomorphisms.
\end{enumerate}
\end{lemma} Thus we have homomorphisms:
$$ \bj :\LLg\to\ho(L,\L^2 (H/L)),\quad
\wj :\LLg\to\ho(L,\Z ) $$

\begin{proof}[Proof of Lemma~\ref{lem.cross}] To prove (a) consider a
generating element $[h_1 ,h_2 ]\in F_2$. If $f\in\LLg$, then we can write
$f_{\ast}(h_i )=h_i l_i$, for some elements $l_i\in\fl$ and so we have:
$$f_{\ast}([h_1  ,h_2 ])=[h_1  ,h_2 ][l_1 ,h_2 ][h_1  ,l_2 ][l_1 ,l_2 ]\
\bmod F_3
$$ Thus $\hj (f)[h_1  ,h_2 ]=l_1\wedge h_2 +h_1 \wedge l_2 +l_1\wedge l_2$.
Since
this element projects to zero in
$\L^2 (H/L)$, (a) follows for
$\bj$. For
$\wj$ we note that $l_1 \cdot l_2 =0$, since the intersection form is zero on
$L$ and $h_1 \cdot l_2 +l_1\cdot h_2 =0$ since $f_{\ast}$ is an isometry.

To prove (b) we need only note that $g_{\ast}$ is the identity on $L$ and
apply (a) and Equation~\eqref{eq.cross}.
\end{proof}

Now we follow Johnson and reformulate $\bj$ first as a homomorphism
$$\tj :\LLg\to H/L\otimes\L^2 (H/L)$$
 using the duality isomorphism $\ho(L,\Z )\iso H/L$ defined by the
intersection form on $H$. If $\{ x_i ,y_i\}$ is an admissible basis of
$F$ for
$L$, then
$\tj (f)=\sum_i y_i\otimes\bj (f)\cdot x_i$.

We have, for any free-abelian group $V$, a short exact sequence:
\begin{displaymath}
\begin{CD} 0 @>>> \L^3 V @>\eta >> V\otimes\L^2 V @>\theta >> \D_3 (V) @>>> 0
\end{CD}
\end{displaymath}
 where $\eta (v_1\w v_2\w v_3 )=v_1\otimes (v_2\w v_3 )+v_2\otimes (v_3\w v_1
)+v_3\otimes (v_1\w v_2 )$ and $\theta (v_1\otimes (v_2\w v_3 )=[v_1 ,[v_2
,v_3 ]]$. We leave the proof as an exercise for the reader.

We can use this sequence to show that  $\tj (\LLg )\sub\L^3
(H/L)$. If
$f\in\LLg$, then we can write $f_{\ast}(x_i )=x_i C_i$, for some $C_i\in F_2$
which represents $\bj (f)\cdot x_i$, and
$f_{\ast}(y_i )=y_i l_i$ for some $l_i\in \fl$. By distributivity of
brackets we have:
\begin{equation}\lbl{eq.dist} f_{\ast}[x_i ,y_i ]\con [x_i ,y_i ][C_i ,y_i
][x_i ,l_i ]\a_i\ \bmod F_4
\end{equation} where $\a_i$ is a product of elements of $\F_3 /\F_4\iso\D_3
(H)$ which vanish when projected into $\D_3 (H/L)$. Since any element of
$\Gg$ leaves the boundary curve of $\Sg$ fixed, we have that, up to
conjugacy,
$f_{\ast}(\prod_i [x_i ,y_i ])=\prod_i [x_i ,y_i ]$. Our first deduction from
Equation~\eqref{eq.dist} is that $\prod_i [x_i ,l_i ]\in F_3$. Now if we
write out
$l_i =(\prod_j x_j^{e_{ij}})C'_i$, for some $C'_i\in F_2$, we can conclude
that
$e_{ij}=e_{ji}$ and so
$\prod _i [x_i ,l_i ]\ (\bmod F_4)$ reduces to a product of triple brackets,
each of which has some $x_j$ as one of its entries, and so vanishes in $\D_3
(H/L)$. Thus we can conclude, from Equation~\eqref{eq.dist}, that $\sum_i
[C_i ,y_i ]=0$ in
$\D_3 (H/L)$. But this means that $\theta\tj (f)=0$ and so $\tj$ actually
gives us a homomorphism $\LLg\to\L^3 (H/L)$.

We can now set $\JJ =\tj\oplus\wj$, where $\wj$ is reformulated as a
homomorphism
$\wj :\LLg\to H/L$, using the duality isomorphism $\ho (L,\Z )\iso H/L$..

\subsection{Proof of (a)}  Since the Johnson homomorphism
$J$ is onto $\L^3 H$, it is obvious that $\tj$ is onto. To deal with
$\wj$ we first point out that $\wj |\Tg$ can be alternatively defined by the
formula: $\wj (f)=\ss J(f)$, where $\ss :\L^3 H\to H\to H/L$ is defined by
contraction:
$$\ss (h_1\w h_2\w h_3  )=(h_1\cdot h_2 )h_3 +(h_2\cdot h_3 )h_1  +(h_3\cdot
h_1  )h_2 $$ Let $\{ x_i ,y_i\}$ be an admissible basis of $H$. Now choose
$f\in\Tg$ so that
$J(f)=x_i\w y_i\w y_j$, for a prescribed $i\not= j$. Then $\wj (f)=y_j$ while
$\tj (f)=0$ in
$\L^3 (H/L)$. From this it follows that $\JJ$ is onto.
\qed

\subsection{Proof of (b)}\lbl{sec.lh}
 It is essentially proved in~\cite{GL} that
$\tj (\Lg )=0$ but we take a different approach to give a unified proof that
$\JJ (\Lg )=0$.

Suppose that $f$ is a Dehn twist along the curve $\g$ in $\Sg$. If $\g$
bounds in
$\Sg$, then $f\in \Kg$ and there is nothing to prove. If $\g$ does not bound,
then we may assume that $\g$ is a meridian curve of a handle in some
representation of $\Sg$ as the boundary of a handle-body. Let $\{ x_i
,y_i\}$ be the standard basis of $F$, where $x_i$ represents the meridian of
the $i$-th handle and $y_i$ represents the longitude. Note that this is not
necessarily an admissible basis. We may assume that
$\g$ represents $x_1$. Then we see easily that
$$ f_{\ast}(x_i )=x_i\ (1\le i\le g),\ f_{\ast}(y_i )=\cases y_i &\text{ if
} 1\le i\le g \\ y_1 x_1 &\text{ if } i=1 \endcases $$ All we know about
$L$, though, is that $x_1\in\fl$ and, as a consequence, an element $w\in F$
lies in $\fl$ only if the exponent sum of $y_1$ is $0$. From this it follows
that for any
$w\in\fl$, $\tj (f)\cdot w=a\w x_1$, where
$a\in H$ is a linear combination of $x_1 ,\cdots ,x_g ,y_2 ,\cdots ,y_g$. But
then it is clear that $a\w x_1$ maps to $0$ in $\L^2 (H/L)$ and
$c(a\w x_1 )=0$. Thus $\Lg\sub\ker\JJ$

Now suppose $f\in\LLg$ and $\JJ (f)=0$. We want to show that $f\in\Lg$. We
first recall that the images of $\LLg$ and of $\Lg$ under the canonical map
$\Gg\to Sp(\Z ,g)$ coincide and are equal to the subgroup $B_g\sub Sp(\Z ,g)$
defined in Section~\ref{sec.L}. This is proved in~\cite{GL}. Since there
exists some $g\in\Lg$ so that $f_{\ast}=g_{\ast}$ on
$H$, then, by replacing $f$ by $fg^{-1}$, we may assume that $f\in\Tg$. Since
$\Kg\sub\Lg$ it will suffice to find $g\in\Lg\cap\Tg$ so that
$J(g)=J(f)$. In fact we will find such a $g\in [\Lg ,\Tg ]\sub\Lg\cap\Tg$
(since
$\Lg$ is invariant under conjugation by $\Tg$), which will then, in addition,
show that
$\Lg\cap\Tg =[\Lg ,\Tg ]\cdot\Kg$ (part of the conclusion of
Theorem~\ref{th.lcL}). Let $K=\ker\{p\oplus c:\L^3 H\to\L^3 (H/L)\oplus
H/L\}$. We need to show that, for any $a\in K$, there exists $g\in [\Lg ,\Tg
]$ such that
$J(g)=a$. Some examples of elements of $K$ are, in terms of an admissible
basis of $H$:
\begin{enumerate}
\item $a=x_i\w x_j\w x_k$ for any $i,j,k$,
\item $a=x_i\w x_j\w y_k$ for any $i,j,k$,
\item $a=x_i\w y_j\w y_k$ for $i,j,k$ distinct and
\item $a=y_i\w x_i\w y_k +x_j\w y_j\w y_k$ for $i,j,k$ distinct.
\end{enumerate} In fact we will leave it as an exercise for the reader to
prove that any element of $K$ is a linear combination of elements of these
four types.

To compute $J$ on $[\Lg ,\Tg ]$ we need the following observation. Let
$f\in\Gg$ and $h\in\Tg$. Then we have, for the classical Johnson homomorphism
$J:\Tg\to\L^3 H$:
\begin{equation}\lbl{eq.jcom} J([f,h])=J(fhf^{-1})-J(h)=(f_{\ast}-1)J(h)
\end{equation} Since $J$ is onto and $f_{\ast}$ can be any element of
$B_g$, for
$f\in\Lg$, this formula determines $J([\Lg ,\Tg ])$. We construct some
examples.
\begin{enumerate}
\item Choose $h$ so that $J(h)=x_i\w x_j\w y_k$ for $i,j$ distinct,  and
$f$ so that $\fs (y_k )=y_k +x_k$ and $\fs$ is the identity on every other
basis element. Then it is straightforward to compute, from
Equation~\eqref{eq.jcom}, that $J([f,h])=x_i\w x_j\w x_k$.
\item Choose $h$ so that $J(h)=y_i\w x_j\w y_k$ for $i,k$ distinct, and
$f$ so that $\fs (y_i )=y_i +x_i$ and $\fs$ is the identity on every other
basis element. Then $J([f,h])=x_i\w x_j\w y_k$.
\item Choose $h$ so that $J(h)=y_i\w y_j\w y_k$, for $i,j,k$ distinct, and
$f$ so that
$\fs (y_i )=y_i +x_i$ and $\fs$ is the identity on every other basis element.
Then $J([f,h])=x_i\w y_j\w y_k$.
\item Choose $h$ so that $J(h)=y_i\w y_j\w y_k$, for $i,j,k$ distinct, and
$f$ so that $\fs (y_i )=y_i +x_j ,\ \fs (y_j )=y_j +x_i$ and$\fs$ is the
identity  on every other basis element. Then
$$  J([f,h])=y_i\w x_i\w y_k +x_j\w y_j\w y_k +x_j\w x_i\w y_k $$
\end{enumerate}

Now it is clear that the examples given here cover all the cases above.

This completes the proof of (b). In fact these arguments also prove (c) and
so the proof of Theorem~\ref{th.L} is complete.
\qed

\subsection{Proof of Corollary~\ref{cor.cl}}
Recall the exact sequence (see e.g.~\cite{M2})
\begin{equation*}
\begin{CD}
1 @>>> \pi_1 (T) @> \eta >> \Gg @> e >> \gg @>>> 1
\end{CD}
\end{equation*}
where $T$ is the tangent circle bundle of $\sg$. The exact homotopy sequence of
the bundle \mbox{$T\to\sg$} gives a central extension (since $\sg$ is
orientable).
$$\begin{CD}
1@>>> \Z @>>> \pi_1 (T) @>\rho >> \pi @>>>1
\end{CD}$$
  where
$\pi =\pi_1 (\sg )$. $\eta$ satisfies the following properties:
\begin{enumerate}
\item $\im\eta\sub\Tg$
\item $\eta\vert\Z$ is defined by
$\eta (1)=$ Dehn twist along a circle bounding the puncture, and $\eta |\pi$ is
characterized by the property that the automorphism of $\pi$ induced by
$\eta (\a)$
is conjugation by $\a$.
\item There is a commutative diagram
$$\begin{CD}
\pi_1 (T) @>\eta >> \Gg \\
@V\rho VV              @VV\xi V \\
\pi @>\eta '>> \aut\pi
\end{CD}$$
where $\aut\pi$ is the group of automorphisms of $\pi$, $\eta ' (g)=$
conjugation by $g$, and $\xi (f)$ is the automorphism of the fundamental
group of
$\sg$ induced by $f$.
\end{enumerate}

\begin{claim}
\begin{enumerate}
\item[\rm{(a)}] $\bj\circ\eta =0$
\item[\rm{(b)}] If $\rho (\a )\in\pi$ has homology class $h\in H$, then
$\wj\eta (\a
)\con h\mod L$.
\end{enumerate}
\end{claim}
Assuming this, we then see, from (a), that $\bj$ induces the homomorphism $J_0$
(since $\llg =e(\LLg )$) and, since, by (b), $\wj\circ\eta$ is onto (and $\lg
=e(\Lg )$), $\ker J_0 =\lg$.

\begin{proof}[Proof of Claim.] It follows from (3) above that $\hj\eta (\a
)\cdot
x=h\w x$, for any $x\in L$, where $h\in H$ is the homology class  of $\rho
(\a )$.
(In fact, since $\im\eta\sub\tg, \hj (\eta (\a ))$ is given by the standard
Johnson
homomorphism and this formula is true for any $x\in H$.) If
$x\in L$, then
$h\w x\to 0\in\L^2 (H/L)$-- this proves (a). Similarly $\wj\eta (\a )\in\ho
(L,\Z )$
is the functional
$x\to h\cdot x$ which, under the duality  $\ho (L,\Z )\iso H/L$,
corresponds to the
reduction of $h$. This proves (b).
\end{proof}

Finally note that $p(H)=0$ since $p(h)$ is defined to be $\sum_i x_i\w
y_i\w h$, for any symplectic basis $\{ x_i ,y_i\}$ of $H$ and, since
we can choose this basis so that $x_i\in L$, each term in this sum
goes to $0\in\L^3 (H/L)$. This completes the proof of Corollary \ref{cor.cl}.
\qed
\subsection{Proof of Corollary~\ref{cor.pl}}
\begin{claim}
Suppose $f'\in\Gg$ induces $f\in\gg$. Then $f'_{\ast}$ (the induced
automorphism of $F$) satisfies
\begin{equation}\lbl{eq.pl}
f'_{\ast}|F_L\con\text{id}\mod [F,F_L ]
\end{equation}
if and only if $f_{\ast}|\pi_L\con\text{id}\mod [\pi ,\pi_L ]$.
\end{claim}
This follows from the fact that $\ker\{ F\to\pi\}\sub [F,F_L ]$.

Now \eqref{eq.pl} is clearly equivalent to $f'\in\LLg$ and $\bj
(f')=0$, and this, by definition, is the same as $J_0 (f)=0$. So
Corollary~\ref{cor.pl} is an immediate consequence of Corollary~\ref{cor.cl}.
\qed
\section{Proof of Theorem~\ref{th.lcL}} As pointed out above, we have already
proved the last assertion of Theorem~\ref{th.lcL} and so we address the first
part. Note that this is equivalent to the assertion that $J((\Lg )_m )=K_m$
for
$m\ge 2$.

\subsection{The case $m=2$} We first show that $J((\Lg )_2 )\sub K_2$. This
is equivalent to the two inclusions:
\begin{enumerate}
\item $\Phi_f (H)\sub\ker\{\L^2 H\to\L^2 (H/L)\}$ and
\item $\Phi_f (L)\sub\L^2 L$
\end{enumerate} where $\Phi_f$ is defined in Section~\ref{sec.L}. (1) is
already proved in~\cite{GL}, so we address (2). According to
Theorem~\ref{th.L}, $\hj$ induces
a function $\Lg\to\ho(\rl ,X)$, where $X=\ker\{\L^2 H\to\L^2 (H/L)\}$. Now
$\L^2 L\sub X$ and the induced function \mbox{$J':\Lg\to\ho(\rl ,X/L^2 L)$}
is a {\em homomorphism}. This follows from an alternative version of the
crossed-homomorphism formula for $\hj$:
$$ \hj (fg)=\hj (f)+\fs\circ\hj (g) $$ where $\fs :\L^2 H\to\L^2 H$ is
defined by $\fs (h_1 \w h_2 )=f_{\ast}(h_1 )\w f_{\ast}(h_2 )$. If
$f\in\LLg$ then $\fs$ induces the identity map on $L$ and on $H/L$. From
this it follows that
$\fs |X\con 1\ \bmod\L^2 L$ and so $J'$ is a homomorphism. But then we can
conclude that $J'((\Lg )_2 )=0$, which means that $\Phi_f (\rl )\sub\L^2 L$
for any $f\in (\Lg )_2$. This proves (2).

We now need to show that  $K_2\sub J((\Lg )_2 )$. In fact we will show that
$K_2\sub J([\Lg\cap\Tg ,\Lg ] )=J([[\Lg ,\Tg ],\Lg ])$. Let $\{ x_i ,y_i\}$
be an admissible basis of $H$. Then $K_2$ is generated by the following
elements:
\begin{itemize}
\item[(a)] $x_i\w x_j\w x_k$, for $i,j,k$ distinct, and
\item[(b)] $x_i\w x_j\w y_k$, for $i\not= j$ and $j\not= k$.
\end{itemize} Choose $f\in\Lg$ so that $\fs (y_i )=y_i +x_i$ and $\fs$ is the
identity on all other basis elements. By Section~\ref{sec.lh} we can choose
$h\in\Lg\cap\Tg$ so that $J(h)=y_i\w x_j\w x_k$, if $i,j,k$ are
distinct. Then, by Equation~\eqref{eq.jcom}, $J([f,h])=x_i\w x_j\w x_k$.
If
$i\not= j$ and $j\not= k$, then, by Section~\ref{sec.lh}, we can choose
$h_2\in\Lg\cap\Tg$ so that $J(h_2 )=y_i\w x_j\w y_k$. Then by
Equation~\eqref{eq.jcom}, $J([f,h_2 ])=x_i\w x_j\w y_k$.

This completes the proof of Theorem~\ref{th.lcL} for the case $m=2$.

\subsection{The cases $m\ge 3$}\lbl{sec.m} We first show that $J((\Lg )_m
)\sub K_m$. We have already proved this for
$m=2$ in the previous section. We now argue by induction on $m$ using
Equation~\eqref{eq.jcom} and the fact that $(\fs -1)K_m\sub K_{m+1}$. This
latter inclusion follows from the expansion:
$$ (\fs -1)(a_1\w a_2\w a_3 )=(\fs -1)a_1\w\fs a_2\w\fs a_3 +a_1\w (\fs
-1)a_2
\w\fs a_3 +a_1\w a_2\w (\fs -1)a_3   $$ and the following easy fact:
\begin{equation}\lbl{eql}
\text{If }  f\in\LLg , \text{ then } (\fs -1)a\in L \text{ for any } a\in H,
\text{ and }(\fs -1)a=0 \text{ if }a\in L
\end{equation}

All that remains to prove is that $K_3\sub J((\Lg )_3 )$. But $K_3 =\L^3 L$
is obviously generated by the elements $x_i\w x_j\w x_k$ and if we choose
$h\in (\Lg )_2$ so that $J(h)=x_i\w x_j\w y_k$, which we can do by the
previous section, and choose $f\in\Lg$ so that $\fs (y_k )=y_k +x_k$ and
$\fs$ is the identity on all other basis elements, then $J([f,h])=x_i\w x_j\w
x_k$

The proof of Theorem~\ref{th.lcL} is now complete.
\qed
\section{Proof of Theorems~\ref{th.lcl} and~\ref{th.lcsL}} The argument is a
generalization of that in Section~\ref{sec.m}. Both these theorems will
follow from:
\begin{lemma}
\begin{enumerate}
\item If $ f\in\Gg$ and $ g\in\Gg [n]$, then $ J_n ([f,g])=(\fs -1)J_n (g)$.
\item If $ f\in\LLg$ , then $ (\fs -1)\F_n^r\sub\F_n^{r+1}$ .
\end{enumerate}
\end{lemma}
\begin{proof}[Proof of (1)] This just follows from the easy formula $ J_n
(fgf^{-1})=\fs J_n (g)$  and the fact that $ J_n$  is a homomorphism.
\end{proof}
\begin{proof}[Proof of (2)] First we show that $ (\fs -1)\D_n^r
(H)\sub\D_n^{r+1}$ . $ \D_n^r$  is generated by $n$-fold brackets
$ [a_1 ,\cdots ,a_i \cdots ,a_n ]$ , where at least $ r$  of the $ a_i$
belong to
$ L$ . We can expand in the following manner:
\begin{equation}\lbl{eq.exp} (\fs -1)[a_1 ,\cdots ,a_i\cdots ,a_n
]=\sum_{i=1}^n [a_1 ,\cdots ,a_{i-1},(\fs -1)a_i ,\fs a_{i+1}\cdots ,\fs a_n
]
\end{equation} Now, it follows from~\ref{eql} that each of the terms on
the right side of~\ref{eq.exp} lies in $ \D_n^{r+1}$.

Now suppose $ a\otimes\a\in\F_n^r$ . Then we can expand
\begin{equation}\lbl{eq.ex} (\fs -1)(a\otimes\a )=(\fs -1)a\otimes\fs\a
+a\otimes (\fs -1)\a
\end{equation}
 If $ a\in L$  and $ \a\in\D_n^{r-1}$ , then the first term on the right
side of~\ref{eq.ex} vanishes and the second lies in $
L\otimes\D_n^r\sub\F_n^{r+1}$ . If
$ \a\in\D_n^r$ , then the first term lies in $
L\otimes\D_n^r\sub\F_n^{r+1}$  and the second lies in $
H\otimes\D_n^{r+1}\sub\F_n^{r+1}$ .
\end{proof} This completes the proof of the Lemma and, therefore, of
Theorems~\ref{th.lcl} and~\ref{th.lcsL}.
\section{Proof of Theorems~\ref{th.br}, \ref{th.br1} and \ref{th.k}}
\subsection{Proof of Theorems~\ref{th.br} and \ref{th.br1}} The framed pure
braid group
$\pg$ on
$g$ strands is defined to be the group, under stacking, of pure braids,
where each strand is equipped with a normal framing which is standard on the
boundary. This is easily seen to be canonically isomorphic to
$\pg^0\times\Z^g$, where $\pg^0$ is the usual pure braid group, since there
is a well-defined self-linking number for a normal framing. the imbedding
$\psi :\pg\to\Gg$ is defined in~\cite{Od} as follows. Choose an imbedding
of $\dg$,
the disk with $g$  holes,
into $\Sg$ so that the internal boundary circles map to $g$ independent
meridians and the external circle maps to a bounding closed curve
(see Figure~\ref{fig1}).

According to the arguments of Artin, $\pg$ can be regarded as the
group of isotopy classes of diffeomorphisms of $\dg$ which are the identity
on the boundary. This representation of $\pg$, together with the chosen
imbedding $\dg\sub\Sg$, defines $\psi$. Of course $\psi$ depends on the
choice of imbedding. An algebraic formulation of $\psi$ can be given as
follows. Let $\{ x_i ,y_i\}$ be a basis of $F=\pi_1 (\Sg )$ which satisfies:
\begin{enumerate}
\item The homology classes of $\{ x_i ,y_i\}$ form a symplectic basis of $H$
\item The boundary circle of $\Sg$ defines the element
$$(x_1\cdots x_g )^{-1}(y_1 x_1 y_1^{-1}\cdots y_g x_g y_g^{-1})\in F$$
\end{enumerate} See Figure~\ref{fig2}.
\begin{figure}[ht]
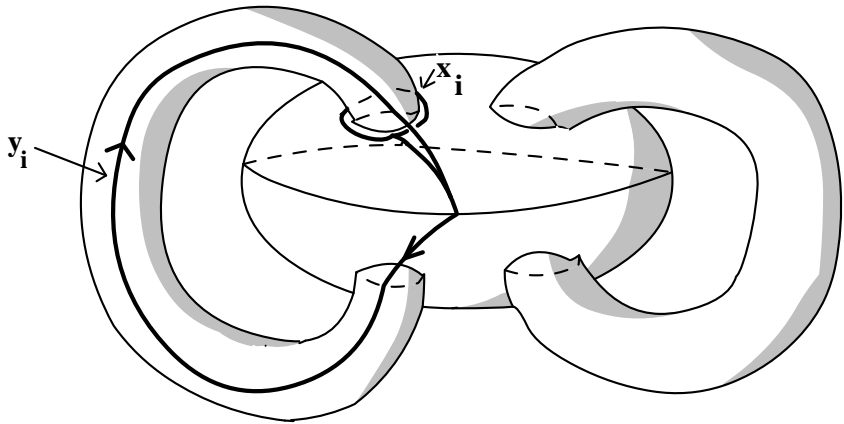

 \TrimLeft{10mm}
 \centerline{\BoxedEPSF{Fig2 scaled 800}}
\bigskip
\caption{Definition of $x_i$ and $y_i$}\lbl{fig2}
\end{figure}
An element $\a\in\pg$ is determined by a sequence $\l_1 ,\cdots ,\l_g\in
F'$, the free group on generators $x_1  ,\cdots ,x_g$, which we shall call
the {\em longitudes } of $\a$. The longitudes satisfy the equation:
\begin{equation}\lbl{eq.long}
\l_1 x_1\l_1^{-1}\cdots\l_g x_g\l_g^{-1}=x_1\cdots x_g
\end{equation} Then $\psi (\a )\in\Gg$ is the element which defines the
following automorphism of $F$:
\begin{equation}\lbl{eq.psi} x_i\to\l_i x_i\l_i^{-1},\quad y_i\to
y_i\l_i^{-1}
\end{equation} Obviously $\psi$ is an imbedding. It is not hard to prove the
following Proposition (see~\cite{O}).
\begin{prop} $\psi$ induces an {\em imbedding }of the quotients:
\begin{displaymath}
\psi_n :(\pg )_n /(\pg )_{n+1}=\pg [n]/\pg [n+1]\to\Gg [n-1]/\Gg [n]
\end{displaymath} for every $n\ge 1$.
\end{prop}
$\{\pg [n]\}$ is the weight filtration of $\pg$, i.e. $\a\in\pg [n]$ if and
only if the longitudes $\l_i\in F'_n$. In fact $(\pg )_n =\pg [n]$.
\begin{proof}[Proof of Theorem~\ref{th.br}.] Since any diffeomorphism of
$\dg$ is a composition of Dehn twists along closed curves in $\dg$, then
$\psi (\pg )\sub\Lg$ if we set $L=\im\{\h (\dg )\to\h (\Sg )\}$. Thus $L$ is
generated by $\{ x_i\}$. Now suppose $\l_i\in F'_n$ are the longitudes of
$\a\in (\pg )_n$. Then
$J_{n-1}\psi_n (\a )\in\ho(H,\D_n (H))$ is given by $x_i\to [\l_i ,x_i
]=0\in\D_n (H),\ y_i\to -\l_i$. The dual element of $H\otimes\D_n (H)$ is,
therefore, $-\sum_i x_i\otimes\l_i$, which belongs to
$L\otimes\D_n (L)$.
\end{proof}
\begin{proof}[Proof of Theorem~\ref{th.br1}.] If $\a\in\pg$ with
longitudes $\l_i\in F'$, then we can write their reductions $l_i\in
L=F'/F'_2$ in the form $l_i =\sum_j a_{ij}x_j$. It is clear that
$a_{ij}$ is the linking number of the $i$-th and $j$-th strands of
(the closure of) $\a$, if $i\not= j$, and the self-linking number of the
$i$-th strand if $i=j$. Thus $A=(a_{ij})$ is a symmetric matrix and it
is clear that any symmetric matrix can be realized by some pure braid. The
definition of $\psi (\a )$ given in~\ref{eq.psi} shows that the
symplectic automorphism determined by $\psi (\a )$ is represented by the
matrix $(\smallmatrix I & -A \\ 0 & I \endsmallmatrix )$. This proves (a).

To prove (b) we will construct an analogue of the Johnson homomorphism $J_b
:(\pg )_2 =\pg [2]\to\L^3 L$ whose kernel is $(\pg )_3 =\pg [3]$ such
that $J_b$ induces an isomorphism $(\pg )_2 /(\pg )_3 \otimes\Q\iso\L^3
L\otimes\Q$. Suppose
$\a\in (\pg )_2$ has longitudes $\l_i\in F'_2$ and, under the isomorphism
$F'_2 /F'_3\iso\L^2 L$, $\l_i\to l_i$. We define $J_b (\a )=\sum_i x_i\otimes
l_i\in L\otimes\L^2 L$. It is clear that the kernel of $J_b$ is $(\pg )_3$.
We see that $\im J_b\sub\L^3 L$ by a copy of Johnson's argument. Consider the
exact sequence:
\begin{equation}
\begin{CD}
0 @>>> \L^3 L @>>> L\otimes\L^2 L @>\sigma >> \D_3 (L) @>>> 0
\end{CD}
\end{equation}
where $\sigma (a\otimes b)=[a,b]$.
Equation~\eqref{eq.long}, when read in $F'_3 /F'_4$, becomes $\prod_i
[x_i ,\l_i ]=1$ or, equivalently, $\sigma J_b (\a )=0$.

Now the desired fact that $\im J_b\otimes\Q =\L^3 L\otimes\Q$ will follow
from the fact that $d_g=\dim ((\pg )_2 /(\pg )_3 \otimes\Q )$ is equal
to $\dim (\L^3 L\otimes\Q )$. To prove this we use the split-exact sequence
which computes the lower central series of the pure braid groups
(see~\cite{F}):
$$ 1\to F'_2 /F'_3 \to (\pg )_2 /(\pg )_3 \to (P_{g-1})_2 /(P_{g-1})_3 \to 1
$$
From this sequence we get the recursive formula $d_g =d_{g-1}+\binom
{g-1}2$. Solving this ($d_1 =0$) gives $d_g =\binom g3$, the dimension
of $\L^3 L$.
\end{proof}
\subsection{Proof of Theorem \ref{th.k}}
First of all it is clear that any simple closed curve in $\dg$ is a
separating curve
in $\Sg$\ , under the imbedding $\dg\sub\Sg$ which defines $\k$. Therefore the
corresponding Dehn twist of $\Sg$ defines an element of $\Kg$. For the
remainder
of the proof we need an algebraic description of the map $\k :\pg\to\Gg$. To
do this we will find it convenient to introduce a map $\d :F'\to F$,  where
$F'$
is the free group on generators$x_1 ,\cdots ,x_g$ and $F$ is the free group on
$x_1 ,\cdots ,x_g ,y_1 ,\cdots ,y_g$, defined by $\d (x_i )=[x_i ,y_i ]$. Now
suppose $\{
\l_i\}$ are the longitudinal elements in $F'$
which define the element $\a\in\pg$. Then $\k (\a )$ is the element of $\Gg$
which corresponds to the automorphism of $F$ given by
\begin{equation}\lbl{eq.k}
\begin{split}
x_i & \to\d (\l_i )x_i\d (\l_i )^{-1} \\
 y_i & \to\d (\l_i )y_i\d (\l_i )^{-1}
\end{split}
\end{equation}
If $\l_i\in F'_n$, then $\d (\l_i )\in F_{2n}$ and so $\k (\a )\in\Gg [2n]$.
But it is proved in \cite{Od} that $\a\in (\pg )_n$ if and only if $\l_i\in
F'_n$ for every $i$. This shows that $\k ((\pg )_n )\sub\Gg [2n]$.

For the proof of injectivity we will need the following:
\begin{lemma}\lbl{lem.comm}
$\d$ induces an injection $F'_n /F'_{n+1}\to F_{2n}/F_{2n+1}$.
\end{lemma}
\begin{proof}  $F'_n /F'_{n+1}$ is the free abelian group generated by the {\em
standard basic Lie elements} $\{ z_{\nu}\}$ on the symbols $x_1 ,\cdots
,x_g$ (with that ordering) of degree $n$, as defined in
\cite[p. 334]{KMS}. It will suffice, therefore, to observe that $\{\d
(z_{\nu})\}$
are distinct standard basic Lie elements on the symbols $x_1 ,\cdots ,x_g ,y_1
,\cdots ,y_g$ (with that ordering) of degree $2n$. But it follows, by a
straightforward inductive argument on the degree that these elements are
basic and,
in the {\em standard ordering}, as defined  in \cite[p. 334]{KMS}, $z_{\nu}<
z_{\mu}$ implies $\d (z_{\nu})< \d (z_{\mu})$.
\end{proof}

To prove that $\k_n$ is injective we use the Johnson-Morita map $J_{2n}:\Gg
[2n]\to
H\otimes\D_{2n+1}(H)$ whose kernel is $\Gg[2n+1]$. If $\a\in (\pg )_n$, then,
from \eqref{eq.k}, we have
$$J_{2n}(\a )=\sum_i (x_i\otimes [\d (\l_i ),y_i ]-y_i\otimes [\d (\l_i
),x_i ]) $$
 It follows from this that if $J_{2n}(\a )=0$, then $\d (\l_i )\in F_{2n+1}$
for all $i$. From Lemma \ref{lem.comm} we conclude that $\l_i\in F'_{n+1}$ and
so $\a\in (\pg )_{n+1}$.

Finally, to see that the induced map $(\pg )_n /(\pg )_{n+1}\to (\Kg )_n /(\Kg
)_{n+1}$ is injective it suffices to observe that $(\Kg )_n\sub\Gg [2n]$, for
every $n$ (see \cite[Cor. 3.3]{M4}), and so the injection $\k_n$ can be
factored
$$ (\pg )_n /(\pg )_{n+1}\to (\Kg )_n /(\Kg )_{n+1}\to\Gg [2n]/\Gg [2n+1] $$
\qed
\section{Proof of Theorem~\ref{th.lkt} and Theorem~\ref{th.ft}}
\subsection{Proof of Theorem~\ref{th.lkt}}
Since $\Kg$ and $\Tg$ are normal and any two choices of $L$ give conjugate
subgroups $\Lg$ the truth of the assertions for any single choice of $L$
implies
the truth for any other choice. Recall from~\cite{GL1} the following material.

Let $i:\sg\to\s3$ be any Heegard imbedding. and $L\sub H$ a suitable Lagrangian
subspace. Then the association $f\to\s3_f$, for any $f\in\LLg$, defines a
linear
map $\Q\LLg\to\M$, where $\Q G$, for any group $G$, denotes the group algebra
over $\Q$. This map induces the following three maps on the associated graded
lower central series quotients:
\begin{equation*}\lbl{eq.phi}
\begin{split}
\phi_n^L &:\G_{3n}(\Lg )\otimes\Q\to\A_n^{\text{conn}}(\emptyset ) \\
\phi_n^T &:\G_{2n}(\Tg )\otimes\Q\to\A_n^{\text{conn}}(\emptyset )\\
\phi_n^K &:\G_{n}(\Kg )\otimes\Q\to\A_n^{\text{conn}}(\emptyset )
\end{split}
\end{equation*}
where $\A_n^{\text{conn}}(\emptyset )$ is the subspace of
$\A_n (\emptyset )$ spanned by {\em connected }trivalent graphs and $\G_n
(G)$ is
the lower central series quotient $G_n /G_{n+1}$, for any group $G$.

It is proved in~\cite{GL1} that these maps are {\em onto} if $g\ge 5n+1$.
Thus we
can, for example, choose $f\in (\Kg )_n$ so that $\phi_n^K (f\otimes 1)\not=
0$. If $f\in (\Lg )_{3n+1}$ then we would have $\phi_n^L (f)$ defined and
zero. But $\phi_n^L$ and $\phi_n^K$ take the same value on any element in both
their domains, since they are defined by the same construction. Thus
$f\notin (\Lg
)_{3n+1}$. Similarly $f\notin (\Tg )_{2n+1}$.

The proofs of the other two non-inclusions are the same.
\qed
\subsection{Proof of Theorem~\ref{th.ft}}
Let $b, D$ and $K\sub\s3$ be as indicated on page \pageref{action}. Let
$N_0\approx I\times D$ be a regular neighborhood of $D$ in $\s3$ and $N$ a
regular neighborhood of $K\cup D$. Let $i:\Sg\to\partial N$ be a
homeomorphism. Then $K_b\sub\s3$ is defined by replacing $K\cap N_0$, which
consists of
$g$ straight line segments, by the pure braid $b$. Thus we have $K_b
\sub N\sub\s3$. Let $\dg =D\times 1\cap\Sg$, which determines the
homomorphism $\Psi$. See Figure~\ref{fig3}.
\begin{figure}[ht]
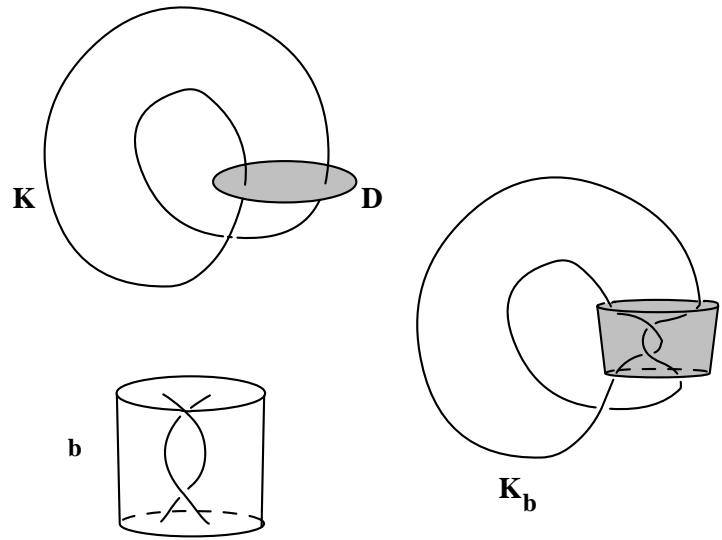

 \centerline{\BoxedEPSF{Fig3 scaled 800}}
\bigskip
\caption{The action of a pure braid on a knot}\lbl{fig3}
\end{figure}

The main point here is the
following:
\begin{lemma}\lbl{lem.act}
There is a homeomorphism $h$ of $N$ onto itself which satisfies:
\begin{enumerate}
\item $h(K)=K_b$,
\item $h|\partial N= i\psi (b)i^{-1}$.
\end{enumerate}
\end{lemma}
\begin{proof}[Proof of Lemma.]
Recall the  correspondence between a pure braid $b$ and a homeomorphism
$h'$ of $D_g$ onto itself. If $b$ is considered as a path in the
configuration space of $g$ points in the $2$-disk, then this path
extends to a diffeotopy of the $2$-disk. The end point of this
diffeotopy fixes $\d D$ and the $g$ points and $h'$ is just the restriction to
$D_g$. The isotopy class, rel $\d D$ of $h'$ is uniquely determined.
For {\em
framed } pure braids, the same construction determines the isotopy class,
rel $\d
D_g$, of
$h'$ $h'$ is the identity on $\d D_g$.
The diffeotopy then defines a diffeomorphism $\bh$ of $I\times D_g$ onto
itself which is the identity on $0\times D_g$, $h'$ on $1\times\dg$
and, if $d$ denotes the union of the $g$ interior disks whose complement
is $\dg$, $\bh (I\times d)=b$. Here we think of a framed pure braid as a
coordinatized tubular neighborhood of the actual braid. We can identify
$N_0$ with $I\times\dg$ so that a tubular neighborhood of $N_0\cap K$
is identified with $I\times d$. Then $\bh$ becomes a diffeomorphism of
$N_0$ and we can extend it over $N$ by declaring it to be the identity
on $N-N_0$. This defines the desired $h$. (1) is clear from the
construction and (2) follows from the definition of $\psi$.
\end{proof}

We now return to the proof of Theorem~\ref{th.ft}.  Let $N'$ be the
closure of the complement of $N$ in $\s3$. Then $\s3_{\psi (b)}$ is, by
definition, $N\cup_f N'$, where $f=i\psi (b)i^{-1}$. But we can use the
diffeomorphism $h$ of Lemma~\ref{lem.act} to construct a diffeomorphism
$\th :\s3_{\psi (b)}\to\s3$ as follows. Regarding $\s3$ as $N\cup_{id}N'$
we define $\th |N=h$ and $\th |N'=$identity. By Lemma~\ref{lem.act}, $\th
(K)=K_b$. The effect of $\th$ on the framing of $K$
is that an $n$-framing on $K$ is mapped to a
$(n+\sum_{ij}l_{ij})$-framing on $K_b$ and so, for our choice of framing
on $b$, a
$+1$-surgery on
$K\sub\s3_{\psi (b)}$ induces, via $\th$, a $+1$-surgery on $K_b\sub\s3$.

Now
$(\s3_{\psi (b)})_K =(\s3_K )_{\psi (b)}$ and when we show that $\psi
(b)\in\L_g^{L_j}$, for the given imbedding  $j:\Sg\sub\s3_K$, the proof of the
Theorem will be complete. Clearly $\im\psi\sub\L_g^{L_i}$ where
$i\Sg\sub\s3$ is
defined from $K$ above. But $L_i\not= L_j$. In fact if $\{ x_i ,y_i\}$ is the
symplectic basis defined by the meridians and longitudes of $\Sg =\d N$--
so that
$L_i$ is spanned by $\{ x_i\}$-- then $L_j$ is spanned by $\{ x_i -y\}$,
where $y=\sum_i y_i$. Thus to show $\psi (b)\in\L_g^{L_j}$ we need to show that
$\psi (b)(y)=y$. But if $\{ l_{ij}\}$ is the linking matrix of the framed
link $b$, then
$\psi (b)(y_i )=y_i +\sum_j l_{ij}x_j $
and so $\psi (b)(y)=y+\sum_j (\sum_i l_{ij})x_j$. Since we have required
$\sum_i
l_{ij}=0$ the proof is complete.
\qed

\ifx\undefined\bysame
	\newcommand{\bysame}{\leavevmode\hbox to3em{\hrulefill}\,}
\fi

\enddocument